\documentclass{article}
\usepackage{amsmath}
\usepackage{amssymb}
\usepackage{amsfonts}
\usepackage[top=2cm,bottom=2cm,left=2.5cm,right=2cm]{geometry}

\usepackage{amsmath}
\usepackage{graphicx}
\usepackage[colorinlistoftodos]{todonotes}
\usepackage[colorlinks=true, allcolors=blue]{hyperref}
\usepackage{amsthm}
\usepackage{bbm}
\usepackage{enumerate}
\usepackage{verbatim}
\usepackage{biblatex}
\usepackage{float}
\usepackage[affil-it]{authblk}

\newtheorem{theorem}{Theorem}[section]
\newtheorem{claim}[theorem]{Claim}

\newtheorem{lemma}[theorem]{Lemma}

\newtheorem{corollary}[theorem]{Corollary}
\newtheorem{definition}[theorem]{Definition}

\bibliography{sector}

\title{Sectorial equidistribution of the roots of $x^2 + 1$ modulo primes}
\author[1]{Evgeny Musicantov\thanks{musicantov.evgeny@gmail.com}}
\author[2]{Sa'ar Zehavi\thanks{saarzehavi@gmail.com}}
\affil[1]{The Hebrew University of Jerusalem, Jerusalem, Israel}
\affil[2]{Tel Aviv University, Tel Aviv, Israel}
\begin{document}

\maketitle

\begin{abstract}
    The equation $x^2 + 1 = 0\mod p$ has solutions whenever $p = 2$ or $4n + 1$. A famous theorem of Fermat says that these primes are exactly the ones that can be described as a sum of two squares. That the roots of the former equation are equidistributed is a beautiful theorem of Duke, Friedlander and Iwaniec from 1995. We show that a subsequence of the roots of the equation remains equidistributed even when one adds a restriction on the primes which has to do with the angle in the plane formed by their corresponding representation as a sum of squares. 
    
    Similar to Duke, Friedlander and Iwaniec, we reduce the problem to the study of certain Poincare series, however, while their Poincare series were functions on an arithmetic quotient of the upper half plane, our Poincare series are functions on arithmetic quotients of $SL_2(\mathbb{R})$, as they have a nontrivial dependence on their Iwasawa $\theta$-coordinate. Spectral analysis on these higher dimensional varieties involves the nonspherical spectrum, which posed a few new challenges. A couple of notable ones were that of obtaining pointwise bounds for nonspherical Eisenstein series and utilizing a non-spherical analogue of the Selberg inversion formula.
\end{abstract}

\section{Introduction}
Let $p$ be a prime number and consider the modular equation
\[
X^2 + 1 = 0(p),
\]
Our focus is on the distribution of the roots of this equation as $p$ varies. Consider the set of tuples $(\nu, p)$ where $\nu$ is a root of $X^2 + 1 = 0(p)$ taken from the interval $[0,p]$ and $p$ varies over all primes for which this equation is solvable. From each such tuple we extract an element $\frac{\nu}{p}$ which we name the \textit{normalized root}. Order the set of tuples according to the right dictionary order and denote by $Y$ the corresponding series of normalized roots, i.e.
\[
Y = (\dfrac{1}{2},\dfrac{2}{5},\dfrac{3}{5},\dfrac{5}{13},\dfrac{8}{13},...).
\]
In~\cite{DFI}, Duke, Friedlander and Iwaniec proved that this sequence is equidistributed in the unit interval in the following sense.
\begin{theorem}
Let $0\le x < y \le 1$, then
\[
\lim_{N\rightarrow\infty}\dfrac{\#\{n\le N: Y_n\in(x,y)\}}{N} = y - x.
\]
\end{theorem}
In order to motivate our variant we would like to suggest a simple parametrization of the (root, modulus) pairs. When the modulus is not necessarily a prime, this set could be parametrized by the primitive lattice points on the first quadrant.
\begin{figure}[htp]
    \centering
    \includegraphics[height=6cm]{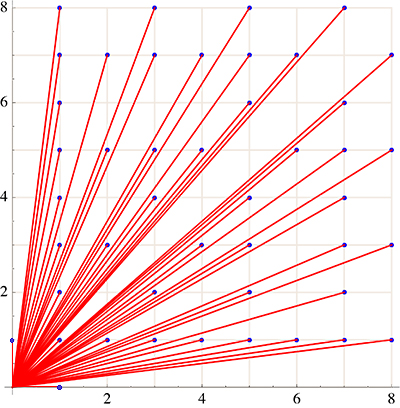}
    \caption{An illustration of a set of primitive lattice points in the first quadrant.}
    \label{fig:primitive}
\end{figure}

For each such primitive lattice point $(a,b)$ on the $\mathbb{Z}_+^2$ plane, denote by $n$ the sum of squares of $a$ and $b$, i.e.
\[
a^2 + b^2 = n.
\]
One may construct a solution to the equation $X^2 + 1 = 0(n)$ by simply considering $\overline{a}b$, where $a$'s inverse is taken relative to multiplication modulo $n$. That $a$ is invertible modulo $n$ follows from the fact that $(a,b)$ is a primitive lattice point. That every (root, modulus)-pair shows up this way exactly once is a simple theorem. That every (root, \textbf{prime} modulus)-pair shows up this way follows from a theorem of Fermat on the decomposition of primes of the form $2$ or $4n+1$ as a sum of two squares.

We call primitive lattice points $(a,b)$ with $a^2 + b^2 = p$ (with $p$ prime) \textit{prime} lattice points. Duke, Friedlander and Iwaniec's theorem says that the sequence of normalized roots corresponding to prime lattice points of radius at most $N$ in \textit{the first quadrant} becomes equidistributed in the unit interval, as $N$ tends to infinity. In this paper, we analyze the sequence of normalized roots to prime moduli, corresponding to prime lattice points taken from a \textit{subsector} of the first quadrant. Our goal is to show that this sequence is also equidistributed in a similar sense.

Several variants of this problem have been studied in the past. For instance, Hooley~\cite{HOOLEY2} demonstrated the equidistribution of the roots of a generic irreducible polynomial of degree greater than one to \textit{composite} moduli. In~\cite{MYARXIV}, the second named author has investigated the joint equidistribution of roots to a pair of polynomial congruences using a similar technique to those in~\cite{HOOLEY2}. Kowalski and Soundararajan~\cite{KOWSOUND} reprove Hooley's theorem (as a special case of a much more general result) in an entirely different fashion, showing that Hooley's theorem is, in fact, an artifact of the Chinese Remainder Theorem.

In~\cite{HOOLEY1}, Hooley uses a different technique to prove a stronger statement regarding the convergence rate (power saving on the relevant Weyl sums), again for composite moduli, but only in the case of an irreducible quadratic polynomial.

\subsection{Our problem}
As mentioned above, the roots of the equations $X^2 + 1 = 0(p)$, with $p$ prime, are parametrized by prime lattice points in the first quadrant. Given two angles $0 \le \alpha<\beta\le \frac{\pi}{2}$, forming a subsector of the first quadrant, we consider a subsequence of $Y$, which consists of these normalized roots to prime moduli for which their corresponding prime lattice point $(a,b)$ lies in the subsector formed by the two angles $\alpha$ and $\beta$ in the first quadrant.

Given such $0 \le \alpha < \beta \le \pi/2$, we obtain a new sequence which we denote by $Y_{\alpha,\beta}$. The goal of this paper is to prove that this sequence is also equidistributed in the unit interval.
\begin{theorem}[Main Theorem]
\label{theorem:main_theorem}
Let $0 \le x < y \le 1$ and $0 \le \alpha < \beta \le \pi/2$. Then
\[
\lim_{N\rightarrow\infty}\dfrac{\#\{n\le N: (Y_{\alpha,\beta})_n\in(x,y)\}}{N} = y - x.
\]
\end{theorem}

\subsubsection{Some numerics}
The following graph shows the percentage of normalized roots to prime moduli, which lie in the interval $[0,0.75]$, taken out of the sequence of normalized roots to prime moduli with prime lattice point in the sector $[0,\pi/6]$.
\begin{figure}[H]
\begin{tabular}{ll}
\includegraphics[width=.4\linewidth]{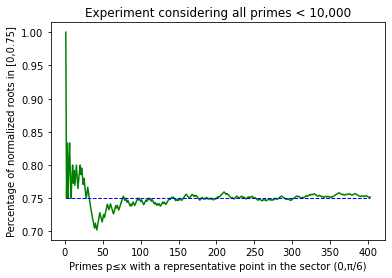}
&
\includegraphics[width=.4\linewidth]{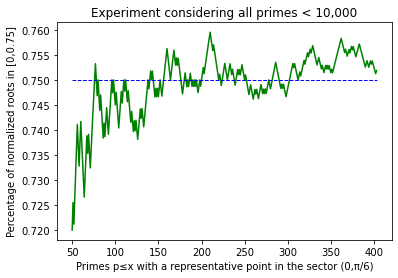}
\end{tabular}

\label{Fig:Race}
\end{figure}
\textbf{Caption}. Left: The green line encodes the percentage of normalized roots for primes $p \le x$ that lie in the interval $[0,0.75]$ out of the $610$ primes that are 2 or 1 mod 4 less than 10,000, under the constraint that their representative prime lattice point lies in the sector $[0,\pi/6]$. There are only 403 such primes. The dashed blue line denotes the conjectured limit. Right: A plot of the graph's tail (from the 50'th data point).

\subsection{The automorphic connection}
\label{section:automorphic_connection}
Gauss observed that, see~\cite{BYKOVSKI}, the (root representative, modulus)-pairs can be parametrized by the orbit of the left action of $SL_2(\mathbb{Z})/\text{Stab}(i)$ on $i$ in the upper half plane, where $\text{Stab}(i)$ is the stabilizer of $i$. Recall that $SL_2(\mathbb{Z})$ acts on the upper half plane by M\"{o}bius transformations, i.e. if $\sigma\in SL_2(\mathbb{Z})$ is such that
\[
\sigma = 
\begin{pmatrix}
a& b\\
c& d
\end{pmatrix},
\]
then it acts on a point $z\in\mathbb{H}$ by
\[
\sigma(z) = \dfrac{az + b}{cz + d}.
\]

The way the root representative and modulus are interpolated from $\sigma$ are by considering the real and imaginary parts of $\sigma(i)$, for example:
\[
\sigma(i) = \dfrac{ai + b}{ci + d} = \dfrac{(ai + b)(-ci + d)}{c^2 + d^2} = \dfrac{ac + bd + i}{c^2 + d^2}.
\]
It may be verified that $ac + bd$ is a root of $X^2 + 1 = 0$ modulo $c^2 + d^2$. This makes $(ac + bd, c^2 + d^2)$ into a (root representative, modulus)-pair, and we regard this pair as the pair corresponding to $\sigma$. Conveniently, a representative of the normalized root is given by the real part of $\sigma(i)$, and the imaginary part of $\sigma(i)$ is the reciprocal of the modulus.

\subsubsection{Moving from root representatives to roots}
Note that in the above parametrization, we made a distinction between ``root representatives" and ``roots". Simply put, given a solution $\nu$ to the equation $X^2 + 1 = 0(n)$, all pairs $(\nu + kn, n)$, with $k\in\mathbb{Z}$, come up from different elements of $SL_2(\mathbb{Z})/\text{Stab}(i)$. However, since all these elements differ by left multiplication by an element of $\Gamma_{\infty}$, replacing the right cosets of $SL_2(\mathbb{Z})/\text{Stab}(i)$ with the double cosets $\Gamma_{\infty}\char`\\SL_2(\mathbb{Z})/\text{Stab}(i)$ turns this into a parametrization of the set of (root, modulus)-pairs.

\subsubsection{Extracting the angle datum}
We explained how to read the (root, modulus) data from an element of the double quotient $\Gamma_{\infty}\char`\\SL_2(\mathbb{Z})/\text{Stab}(i)$, but for our purpose we also require the angular datum corresponding to the associated primitive lattice point. It so happens that this information, including the previous one, sits in the Iwasawa coordinates of the double coset representatives.
\begin{theorem}
There is a bijection between the set of 3-tuples of (normalized root, modulus, angle modulo $\pi/2$) and the set $\Gamma_{\infty}\char`\\SL_2(\mathbb{Z})/\text{Stab}(i)$.
\end{theorem}

The bijection is given by the Iwasawa coordinates of such matrices. The Iwasawa decomposition reduces a matrix $\sigma\in SL_2(\mathbb{Z})$ into 3 coordinates, $(x,y,\theta)$, where the $x$ and $y$ coordinates are equal to ``the normalized root" and ``the reciprocal of the modulus", respectively, while the $\theta$ coordinate corresponds to the angle (of the associated primitive lattice point).

\subsection{Weyl sums and smooth summations}
\label{subsection:weyl_sums}
When trying to prove the equidistribution of a sequence of real numbers $X := (X_1,X_2,...)$ modulo 1, it is often useful to consider the following criterion due to Weyl.
\begin{theorem}[Weyl criterion,~\cite{WEYL}, page 7]
\label{theorem:weyl_criterion}
The sequence $X := (X_1,X_2,...)$ is uniformly distributed modulo 1 if and only if
\[
\lim_{N\rightarrow\infty}\dfrac{1}{N}\sum_{n\le N}e(hX_n) = 0 \; \text{For all integers } h\neq 0,
\]
where $e(z) = e^{2\pi i z}$.
\end{theorem}
\textbf{Remark}: since interchanging $h\leftrightarrow -h$ results in conjugation of the above sum, we may assume $h\in\mathbb{N}$ for simplicity.

Our aim is to prove Theorem~\ref{theorem:main_theorem} using Weyl's criterion. In~\cite{DFI}, Duke, Friedlander and Iwaniec show how the conditions of Weyl's criterion for a sequence of normalized roots to prime moduli can be inferred from certain uniform bounds on Weyl sums of normalized roots whose moduli lie in arithmetic progressions. Such Weyl sums in arithmetic progressions are called ``linear forms".

We make the following definition.
\begin{definition}[Linear forms]
\label{definition:linear_sums}
Let $h,q\in\mathbb{N}$. We define the Weyl sum of parameter $h$ and arithmetic progression of difference $q$ to be
\[
\mathcal{L}_{h,q}(N) = \sum_{qn\le N}\rho_h(qn),
\]
where $\rho_h(n) = \sum_i e(\dfrac{\nu_i}{n})$, and the sum is on those normalized roots $\left(\dfrac{\nu_i}{n}\right)\in Y_{\alpha,\beta}$ with modulus $n$.
\end{definition}

\textbf{Remark}: As in~\cite{DFI}, we are also interested in bilinear sums, but omit their introduction from this exposition.

\subsubsection{Duke, Friedlander and Iwaniec's idea}
Following the brilliant idea of Duke, Friedlander and Iwaniec, we would like to analyze the linear forms $\mathcal{L}_{h,q}(N)$. We start with the case of $q = 1$, and denote by $\mathcal{L}^{DFI}_{h,q}(N)$ the linear forms for the sequence $Y$ of normalized roots without any constraints on the angle corresponding to the root.

Based on the discussion in the previous section regarding the parametrization of the normalized roots using elements of the double coset $\Gamma_{\infty}\char`\\SL_2(\mathbb{Z})/\text{Stab}(i)$, we have
\[
\mathcal{L}^{DFI}_{h,1}(N) = \sum_{\substack{\sigma\in\Gamma_{\infty}\char`\\SL_2(\mathbb{Z})/\text{Stab}(i)\\\frac{1}{y(\sigma)} \le N}}e(hx(\sigma)).
\]
Duke et al. consider the following, smooth version of $\mathcal{L}^{DFI}_{h,1}(N)$
\[
\dfrac{1}{4}\sum_{\sigma\in\Gamma_{\infty}\char`\\SL_2(\mathbb{Z})}e(hx(\sigma))F(4\pi hy(\sigma)),
\]
where upon removing the right quotient by $\text{Stab}(i)$ we now traverse each (root, modulus)-pair with multiplicity $|\text{Stab}(i)| = 4$, thus the factor of $1/4$. 

\textbf{Remark}: here the function $F:\mathbb{R}_+\longrightarrow\mathbb{R}_+$ should be thought of as an indicator over the interval of $y$-values $[1/N,1]$. In practice, Duke, Friedlander and Iwaniec consider a dyadic subdivision of the above smooth summation, as do we.

\begin{definition}[The indicator function $F$]
\label{definition:function_F}
We define the function $F:\mathbb{R}_+\longrightarrow\mathbb{R}_+$ to be any smooth compactly supported function in the interval $[Y/2,Y]$, where $Y = 4\pi h/N$, satisfying the property that for all natural numbers $j \le 12$, $|F^{(j)}| \ll Y^{-j}$.
\end{definition}

\textbf{Remark}: the construction of a smooth indicator function satisfying these properties appears in~\cite[Appendix A]{IWANIECZETA}.

As a next step, we discuss Duke, Friedlander and Iwaniec's representation of the smooth linear forms for a generic difference $q$.

It so happens that for all $\sigma\in \Gamma_0(q)$ and $\tau\in SL_2(\mathbb{Z})$, $\tau$ and $\sigma\tau$ correspond to (root,modulus)-pairs whose modulus have the same value modulu $q$. This motivates a smooth version of the Weyl sums in arithmetic progression of difference $q$, $\mathcal{L}^{DFI}_{h,q}(N)$, via:
\[
\dfrac{1}{4}\sum_{\substack{\tau\in \Gamma_0(q)\char`\\SL_2(\mathbb{Z})\\ y(\tau)^{-1} = 0(q)}}\sum_{\sigma\in\Gamma_{\infty}\char`\\\Gamma_0(q)}e(hx(\sigma\tau))F(4\pi hy(\sigma\tau)).
\]

\subsubsection{Our setting}
\label{subsubsection:our_variant}
In our setting, the constraint on the angle is translated into a constraint on the Iwasawa $\theta$-coordinate. We write the smooth linear forms associated to our problem, denoted $\mathcal{L}_{h,q}^*(N)$, via
\[
\mathcal{L}_{h,q}^*(N) = \dfrac{1}{4}\sum_{\substack{\tau\in \Gamma_0(q)\char`\\SL_2(\mathbb{Z})\\ y(\tau)^{-1} = 0(q)}}\sum_{\sigma\in\Gamma_{\infty}\char`\\\Gamma_0(q)}e(hx(\sigma\tau))F(4\pi hy(\sigma\tau))G(\theta(\sigma\tau)).
\]
where $G(\cdot)$ is to be thought of as an indicator of the interval $[\alpha,\beta]$, defined in a similar way to $F$.

\begin{definition}
\label{definition:function_G}
Given parameters $\alpha,\beta\in \mathbb{S}^1$ as above, and $N^{-1/3} < Z < (\beta - \alpha)/2$, we define the function $G:\mathbb{S}^1\longrightarrow\mathbb{R}_+$ to be any smooth compactly supported function in the interval $[\alpha,\beta]$, satisfying the following properties.
\begin{itemize}
    \item $G$ attains the value $1$ in $[\alpha + Z,\beta - Z]$.
    \item For all $j\in\mathbb{N}$, $||G^{(j)}||^2 \ll_j Z^{1-2j}$.
    \item $G$'s Fourier coefficients, denoted $g_n$, satisfy $|g_n| \ll e^{-2\sqrt{\pi|n|Z}}$.
\end{itemize}
\end{definition}

The reasoning behind the constraint $Z < (\beta - \alpha)/2$ is clear, the constraint $N^{-1/3} < Z$ becomes a necessity from the proof of Claim~\ref{claim:reformulation_large_parameter}. Ultimately, we fix a precise value for $Z$ (in the range $N^{-1/3} < Z < (\beta - \alpha)/2$) in section~\ref{section:equidistribution_to_primes}, through optimization.

The idea is that the above sum can be thought of as a sum over special values of a function on $\Gamma_0(q)\char`\\SL_2(\mathbb{R})$. Let $g\in SL_2(\mathbb{R})$. Denote by $P(g)$ the function defined as
\[
P(g) := \sum_{\sigma\in\Gamma_{\infty}\char`\\\Gamma_0(q)}e(hx(\sigma g))F(4\pi hy(\sigma g))G(\theta(\sigma g)).
\]
In terms of $P(g)$, $\mathcal{L}_{h,q}^*(N)$ can be written as
\[
\mathcal{L}_{h,q}^*(N) = \dfrac{1}{4}\sum_{\substack{\tau\in \Gamma_0(q)\char`\\SL_2(\mathbb{Z})\\ y(\tau)^{-1} = 0(q)}}P(\tau).
\]
The function $P(\cdot)$ is sometimes referred to as a Poincare series, and the majority of this body of work is dedicated to the construction of pointwise bounds for this Poincare series.

\subsection{Rundown of our technique}
\label{subsection:rundown}
As mentioned above, our main theorem essentially follows from pointwise bounds on our Poincare series. Different than the case analyzed in~\cite{DFI}, our Poincare series is not an automorphic function on $\Gamma_0(q)\char`\\\mathbb{H}$, but rather on $\Gamma_0(q)\char`\\SL_2(\mathbb{R})$. In~\cite{DFI}, the authors study their Poincare series by considering its spectral decomposition into Laplace eigenforms. In our case, the Laplacian is replaced by the Casimir, which has a bigger spectrum, and this fact poses some of the challenges which we overcome in this paper.

In section~\ref{section:plancharel}, we explain the structure of the spectrum of the Riemann surface $\Gamma_0(q)\char`\\SL_2(\mathbb{R})$, and introduce the Plancharel formula.

In section~\ref{section:point-wise_bound_overview}, we overview our method for constructing pointwise upper bounds on our Poincare series. Roughly speaking, as a first step, we apply a weighted version of Cauchy-Schwarz, which yields a bound of the form
\[
|P(g)|^2 \ll K(g)\cdot R(P),
\]
where $R(P)$ is a weighted sum of squares of the magnitudes of the projections of the Poincare series $P$ on Casimir eigenfunctions of our modular variety, and $K(g)$ is a weighted sum of squares of the absolute values of the pointwise evaluations of the eigenfunctions of the Casimir on our modular variety at the point $g$. 

In section~\ref{section:the_sum_k_q} we bound the sum $K(g)$ for the values of $g$ we are interested in. Following Duke, Friedlander and Iwaniec, we reduce the sum over the spectrum of the Casimir on $\Gamma_0(q)\char`\\SL_2(\mathbb{R})$ to an analogous sum over the spectrum of the Casimir on $\Gamma_0(1)\char`\\SL_2(\mathbb{R})$. In~\cite{DFI}, a spherical analogue of what we do is obtained using the inverse Selberg transform. This formula however is only applicable in the spherical case, and so we appeal to a non-spherical variant of the inverse Selberg transform, introduced by Hejhal~\cite[p. 386]{HEJHAL2}. Moreover, in~\cite{DFI}, the corresponding trace sum $K(g)$ for the full modular group is bounded naively term by term using uniform bounds in the spectral aspect, and then applying the Selberg trace formula. To get a reduction to the Selberg trace formula, we require bounds that are uniform in both the spectral and weight aspects of our Maass forms. In the cuspidal case, our estimates follow from a result of Bernstein and Reznikov, see~\cite{BERREZ}. In~\cite{MUSICZEHAVI_EISENSTEIN}, we extend Bernstein and Reznikov's techniques for Eisenstein series.

In section~\ref{section:bounding_r(p)}, we bound the sum $R(P)$ by reducing it to a (weighted) sum of squares of absolute values of Fourier-Whittaker coefficients of Maass forms on $\Gamma_0(q)\char`\\SL_2(\mathbb{R})$. This is handled by the Kuznetsov/Petersson trace formulas.

In section~\ref{section:equidistribution_to_primes}, we turn our bound on the smooth Weyl sums into a bound on our (non-smooth) Weyl sums by studying the difference between the two. This reduces to an elementary lattice point counting problem. We then obtain a bound on our Weyl sums to prime moduli using a recipe of T\'{o}th~\cite{TOTH}, which uses the sieving argument taken from~\cite{DFI}, and results in the proof of our main theorem.

\subsection*{Acknowledgements}
This research was supported by the European Research Council (ERC) under the European Union's Horizon 2020 research and innovation programme (Grant agreement No. 786758)

The authors would like to thank Ze\'{e}v Rudnick for his invaluable support, suggestions and very many discussions on this problem. The authors would like to thank Valentin Blomer for taking interest in this project and sharing some helpful comments and references, as well as to Bingrong Huang for many discussions on the theory of automorphic forms.


\section{The Plancherel formula on $X_0(q)$}
\label{section:plancharel}
Denote by $V(X_0(q))$ the Hilbert space of smooth square integrable functions on $X_0(q) = \Gamma_0(q)\char`\\SL_2(\mathbb{R})$ relative to the standard Haar measure ($d\chi = dxdyd\theta/y^2$). This space enjoys a spectral decomposition into eigenspaces of the Casimir, denoted $\mathfrak{C}$, which is the second order differential operator, given in Iwasawa coordinates by
\[
\mathfrak{C} := -y^2\left(\left(\dfrac{\partial}{\partial x}\right)^2 + \left(\dfrac{\partial}{\partial y}\right)^2\right) + y\dfrac{\partial}{\partial x}\dfrac{\partial}{\partial\theta}.
\]
For more information, see Bump~\cite{BUMP}.

We have the Plancherel Formula,
\begin{theorem}[The Plancherel Formula]
\label{theorem:plancherel}
$\forall f\in V(X_0(q)):$
\[
f(g) = \sum_j\sum_n<f,u_{\lambda_j,n}>u_{\lambda_j,n}(g) + \dfrac{1}{4\pi}\sum_{\mathfrak{a}}\int_{-\infty}^{\infty}\sum_{n}<f,E_{\mathfrak{a},n}(*,1/2 + it)>E_{\mathfrak{a},n}(g,1/2 + it)dt.
\]
\end{theorem}
We begin with a description of this formula. The two parts
\[
\sum_j\sum_n<f,u_{\lambda_j,n}>u_{\lambda_j,n}(g)\text{ and }\dfrac{1}{4\pi}\sum_{\mathfrak{a}}\int_{-\infty}^{\infty}\sum_{n}<f,E_{\mathfrak{a},n}(*,1/2 + it)>E_{\mathfrak{a},n}(g,1/2 + it)dt,
\]
correspond to the discrete and continuous spectrum of the Casimir on $X_0(q)$, respectively.

\subsection{The discrete part}
The first double summation,
\[
\sum_j\sum_n<f,u_{\lambda_j,n}>u_{\lambda_j,n}(g),
\] 
consists of the contribution of the discrete spectrum. The outer summation $\sum_{j}$, is a summation on ``certain" Casimir eigenspaces $V_{\lambda_j}$, indexed by their eigenvalues $\lambda_j$ (i.e. $\forall f\in V_{\lambda_j}: \mathfrak{C} f = \lambda_j f$), such that
\[
0 = \lambda_0 < \lambda_1 \le \lambda_2 \le \lambda_3 \le ...
\]
The $V_{\lambda_j}$ Casimir eigenspaces appearing in our decomposition are generally not the entire $\lambda_j$-eigenspace, except for $\lambda_0 = 0$ which consists of just one eigenvector, which is 1. To be precise, the $V_{\lambda_j}$ are the irreducible cuspidal representations, relative to the right action of $SL_2(\mathbb{R})$.

\begin{definition}[The right action of $SL_2(\mathbb{R})$ on $V(X_0(q))$]
The right action of $SL_2(\mathbb{R})$ on $V(X_0(q))$, denoted by $\rho$, is the linear operator $\rho:SL_2(\mathbb{R})\longrightarrow GL(V(X_0(q))$, defined for $g\in SL_2(\mathbb{R})$ and $f\in V(X_0(q))$, by 
\[
(\rho(g)f)(h) = f(hg)
\]
for all $h\in SL_2(\mathbb{R})$.
\end{definition}

In order to describe the inner summation $\sum_n$, we require the definition of ``pure weight" functions.

\begin{definition}[Functions of pure weight]
\label{definition:pure_weight}
A function $f\in C^{\infty}(SL_2(\mathbb{R}))$ is of pure weight $n\in \mathbb{Z}$, if for any $k = k(\theta) \in K = SO_2(\mathbb{R})\subset SL_2(\mathbb{R})$, i.e. $k$ is a $\theta$-radian rotation matrix, one has
\[
\rho(k)f = e^{in\theta}f.
\]
\end{definition}
In other words, pure weight functions are those on which the right translations by rotation matrices act by a (fixed) character.

It is a basic fact that each $V_{\lambda_j}$ is equal to a direct sum $\bigoplus_{n}V_{\lambda_j,n}$, where each $V_{\lambda_j, n}$ is a pure weight $n$ subspace of $V_{\lambda_j}$. It is well known that $\forall j, n\in\mathbb{Z}: dim(V_{\lambda_j,n}) \le 1$, see~\cite{BUMP}. We fix $u_{\lambda_j,n}$ to be any unit vector in $V_{\lambda_j,n}$.

\subsection{The continuous part}
The second term,
\[
\dfrac{1}{4\pi}\sum_{\mathfrak{a}}\int_{-\infty}^{\infty}\sum_{n}<f,E_{\mathfrak{a},n}(*,1/2 + it)>E_{\mathfrak{a},n}(g,1/2 + it)dt,
\]
consists of the contribution of the continuous spectrum. The functions $E_{\mathfrak{a},n}(g,1/2 + it)$ are the Eisenstein series. Denote by $\mathfrak{C}$ the set of cusps of $X_0(q)$. The first summation $\sum_{\mathfrak{a}}$ is a summation over $\mathfrak{a}\in \mathfrak{C}$. The integration is taken over the type parameter $t$, which is related to the corresponding Casimir eigenvalue via $\lambda = 1/4 + t^2$. Together, these two parameters (cusp and type) parametrize the irreducible continuous representations, with respect to right translation by $SL_2(\mathbb{R})$ (up to parity, which we discuss next). The innermost summation is over pure weight $n\in\mathbb{Z}$ Eisenstein series; the space of Eisenstein series of a fixed cusp $\mathfrak{a}$, type $t\in\mathbb{R}$ and weight $n\in\mathbb{Z}$ is exactly 1-dimensional.

Since the Eisenstein series are not square integrable, we cannot normalize our basis functions by their norm, as we do in the case of the discrete spectrum. Instead, we pick our basis vectors such that their constant term (in the Fourier-Whittaker decomposition) is equal to 1.

\subsection{The parity of an irreducible representation}
Another important fact regarding the irreducible representations, both discrete and continuous, is that each irreducible representation has a parity. Denote by $V$ an irreducible representation, and by $V_n$ its (at most 1 dimensional) pure weight $n$ subspace. Then either $dim(V_{n}) = 0$ for all even $n$ or $dim(V_{n}) = 0$ for all odd $n$. In fact, each irreducible representation would have a so called minimal $K$-type, which is the smallest weight $n$ in absolute value, for which $dim(V_{n}) = 1$. In this case both $V_{n}$ and $V_{-n}$ would be one dimensional, and their parity decides the parity of the representation.

\textbf{Remark}: it follows that for each cusp $\mathfrak{a}$ and type $t$, there are 2 irreducible Eisenstein representations, one for even weights and on for odd weights.


\section{Constructing an upper bound on $|P(\tau)|$ - An overview}
\label{section:point-wise_bound_overview}
Our goal in this section is to outline our strategy for bounding $|P(\tau)|$. Recall the definition of our Poincare series,
\[
P(g) = \sum_{\gamma\in \Gamma_{\infty}\char`\\\Gamma_0(q)}e(x(\gamma g))F(4\pi h y(\gamma g))G(\theta(\gamma g)).
\]

\begin{theorem}
\label{theorem:point-wise_bound}
For all $\tau\in SL_2(\mathbb{Z})$, $0 < \delta < 1/4$ and $h < N^{1/3}$ (where $Y = 4\pi h/N$). One has
\[
|P(\tau)| \ll_{\delta} hZ^{-19}Y^{-1-2\delta}\left(1 + h^{\frac{1}{2}}|\log Y|^2\left(Y + Y^{-1}\right)^{\frac{1}{2}}q^{-1}(h,q)^{\frac{1}{2}}\tau(hq)\right).
\]
\end{theorem}
Since $P(g)$ is $\Gamma_0(q)$-left shift invariant, our Poincare series is a function on the modular curve $X_0(q) = \Gamma_0(q)\char`\\SL_2(\mathbb{R})$. Moreover, since it is smooth and compactly supported, it is also square integrable. By the Plancherel Formula~\ref{theorem:plancherel}, we have
\begin{theorem}
\label{theorem:poincare_series_decomposition}
\[
P(g) = \sum_j\sum_n<P,u_{\lambda_j,n}>u_{\lambda_j,n}(g) + \dfrac{1}{4\pi}\sum_{\mathfrak{a}}\int_{-\infty}^{\infty}\sum_{n}<P,E_{\mathfrak{a},n}(*,1/2 + it)>E_{\mathfrak{a},n}(g,1/2 + it)dt.
\]
\end{theorem}
\textbf{Remark}:
this decomposition is generally an $L^2$-equivalence, however, because our Poincare series is smooth and compactly supported, the right hand side is a continuous function. Since a continuous representative of an $L^2$-equivalence is unique (if it exists), the spectral decomposition converges pointwise to $P(g)$.

\begin{claim}
$P(g)$ is an even function.
\end{claim}
\begin{proof}
Because $G(\theta)$ is $\pi/2$-radian rotation invariant, it is also $\pi$-radian rotation invariant, making it an even function.
\end{proof}
\textbf{Remark}: being an even function, means that our Poincare series' spectral decomposition is supported on even weights.

A weighted Cauchy-Schwarz gives
\[
|P(g)|^2 \le K_q(g)\cdot R(P),
\]
where
\[
K_q(g) = \sum_j\sum_n \rho(t_j,n)|u_{\lambda_j,n}(g)|^2 + \dfrac{1}{4\pi}\sum_{\mathfrak{a}}\int_{-\infty}^{\infty}\sum_n \rho(t,n)|E_{\mathfrak{a},n}(g,1/2 + it)|^2dt,
\]
and
\[
R(P) = \sum_j\sum_n\rho(t_j,n)^{-1}|<P,u_{\lambda_j,n}>|^2 + \dfrac{1}{4\pi}\sum_{\mathfrak{a}}\int_{-\infty}^{\infty}\sum_n\rho(t,n)^{-1}|<P,E_{\mathfrak{a},0}(*,1/2 + it)d>|^2dt.
\]
\textbf{Remark}: the term $K_q(g)$ is independent of the function $P$, and the term $R(P)$ is independent of the point $g$.

The positive weight function $\rho(t,n)$ is introduced in order to allow $K_q(g)$ to converge. Since our precise choice of $\rho(t,n)$ is technical and not very illuminating at this point, we postpone its definition for the next section.

\section{Bounding the sum $K_q(\tau)$}
\label{section:the_sum_k_q}
In this section, we prove:
\begin{theorem}
\label{theorem:bound_on_K}
For all $\tau\in SL_2(\mathbb{Z})$, one has
\[
K_q(\tau) \ll 1.
\]
\end{theorem}

\subsection{Preliminary discussion}
\label{section:preliminary_discussion_k}
We require a few definitions, see Hejhal~\cite[p. 357-359]{HEJHAL1}.

The projection $\pi:SL_2(\mathbb{R})\longrightarrow\mathbb{H}$ is given by
\[
\pi(g) = x(g) + iy(g).
\]
Let $g,h\in SL_2(\mathbb{R})$. The hyperbolic distance between $g$ and $h$, denoted $u(g,h)$, is given by
\[
u(g,h) = \dfrac{|\pi(g) - \pi(h)|^2}{4y(g)y(h)}.
\]
Let $z\in\mathbb{H}$, and $\sigma\in SL_2(\mathbb{R})$, given by $\sigma =
\begin{pmatrix}
a& b\\
c& d
\end{pmatrix}$, then the function $j_{\sigma}:\mathbb{H}\longrightarrow\mathbb{C}$ is
\[
j_{\sigma}(z) := \dfrac{cz + d}{|cz + d|} = e^{i\arg (cz+d)}.
\]
Let $n\in \mathbb{Z}$, we define $H_n:\mathbb{H}\times \mathbb{H}\longrightarrow\mathbb{C}$ by
\[
\forall z,w\in\mathbb{H}:\quad H_n(z,w) = i^n\dfrac{(w-\overline{z})^n}{|w-\overline{z}|^n}.
\]
Let $\Phi_n:\mathbb{R}_+\longrightarrow\mathbb{R}$ be any $C^2$-class function. We define the weight-$n$ point-pair invariant attached to $\Phi_n$ to be the function
\[
k_n(g,h) = \Phi_n(u(\pi(g),\pi(h)))H_n(\pi(g),\pi(h))e^{in(\theta_g - \theta_h)}.
\]
\begin{lemma}
\label{lemma:invariancy_lemma}
A (weight-$n$) point-pair invariant has the property that for all $\sigma\in SL_2(\mathbb{R})$,
\[
k_n(\sigma g,\sigma h) = k_n(g,h).
\]
\end{lemma}
The weight $n$ Laplacian on $SL_2(\mathbb{R})$, denoted $\Delta_n$ is given by
\[
\Delta_n = -y^2(\dfrac{d^2}{dx^2} + \dfrac{d^2}{dy^2}) + iny\dfrac{d}{dx}.
\]
The weight-$n$ point-pair invariant has the property:
\begin{claim}[cf.~\cite{HEJHAL1} page 359, Definition 2.10]
\label{claim:scalar_emition}
Let $f\in C^2(SL_2(\mathbb{R}))$ be a pure weight $n$ Laplace eigenfunction, i.e. $f$ satisfies $\Delta_n f = \lambda f = (1/4 + t^2)f$, then
\[
\int_{SL_2(\mathbb{R})}k_n(g,h)f(h)d\chi(h) = \rho_n(t)f,
\]
where $\rho_n(t)$ depends only on $\Phi_n,n,t$. The value of $\rho_n(t)$ is therefore independent of $f$.
\end{claim}
A point-pair invariant $k:SL_2(\mathbb{R})\times SL_2(\mathbb{R})\longrightarrow \mathbb{C}$ is any series
\[
k(g,h) = \sum_{n\in\mathbb{Z}}k_n(g,h),
\]
where the series is absolutely convergent for all $g$ and $h$, and the $k_n(g,h)$ are weight-$n$ point-pair invariants.

A point-pair invariant $k(g,h)$ induces an automorphic kernel $K_q(g,h):SL_2(\mathbb{R})^2\longrightarrow \mathbb{C}$, given by
\[
K_q(g,h) = \sum_{\gamma \in \Gamma_0(q)}k(g,\gamma h).
\]
\textbf{Remark}: since the point-pair invariants we are interested in would induce automorphic kernels with the property $|K_q(g,h)| < \infty$, we assume our automorphic kernel has this property for the rest of this discussion.

The automorphic kernel $K_q(g,h)$ acts linearly by an integral transformation on $C^{2}(X_0(q))$, via
\[
\forall f\in C^{2}(X_0(q)):\quad \int_{X_0(q)}K_q(g,h)f(h)d\chi(h) = \int_{SL_2(\mathbb{R})}k(g,h)f(h)d\chi(h).
\]
The equality follows from the unfolding method, valid whenever one has absolute convergence, i.e. when
\[
\int_{SL_2(\mathbb{R})}|k(g,h)f(h)|d\chi(h) < \infty.
\]

The functions $f\in C^{2}(X_0(q))$ we are interested in are Casimir eigenfunctions on the modular curve $X_0(q)$, which have bounded polynomial growth, and our point-pair invariant $k(g,h)$ will be constructed in a way which forces such absolute convergence. In fact, we will demand it to satisfy the stronger property, that for any Casimir eigenfunction $f\in C^{\infty}(X_0(q))$, one has
\[
\int_{SL_2(\mathbb{R})}\sum_{n\in\mathbb{Z}}|k_n(g,h)f(h)|d\chi(h) < \infty.
\]
This implies that for all pure weight-$n$ Laplace eigenfunctions $f$ on $X_0(q)$, we have
\[
\int_{X_0(q)}K_q(g,h)f(h)d\chi(h) = \rho_n(t)f(g).
\]
From this point onward we assume this absolute convergence property of $k(g,h)$.

Denote $\rho(t,n) = \rho_n(t)$. The previous observation implies that $K_q(g,h)$ has a spectral decomposition of the form
\[
K_q(g,h) =
\sum_j\sum_n \rho(t_j,n)u_{\lambda_j,n}(g)\overline{u_{\lambda_j,n}(h)} + \dfrac{1}{4\pi}\sum_{\mathfrak{a}}\int_{-\infty}^{\infty}\sum_n \rho(t,n)E_{\mathfrak{a},n}(g,1/2 + it)\overline{E_{\mathfrak{a},n}(h,1/2 + it)}dt.
\]
The function $\rho(t,n)$ is also known as the Harish-Chandra/Selberg transform of the point-pair $k(g,h)$.

By Claim~\ref{claim:scalar_emition}, and our absolute convergence assumptions, the scalar $\rho(t,n)$ can be computed directly from the point-pair invariant $k(g,h)$ as the unique complex number satisfying
\[
\int_{SL_2(\mathbb{R})}k(g,h)\phi_{\lambda,n}(h)d\chi(h) = \rho(t,n)\phi_{\lambda,n}(g),
\]
for $\phi_{\lambda,n}$ any pure-weight $n$ Casimir eigenfunction of eigenvalue $\lambda = 1/4 + t^2$. Specializing to the case where $\phi_{\lambda,n}(x,y,\theta) = y^{1/2 + it}e^{in\theta}$, where $\lambda = \lambda_t = 1/4 + t^2$, and fixing $g = I\in SL_2(\mathbb{R})$, the identity element, $\rho_k(t,n)$ is computed directly via:

Given a point-pair invariant $k(g,h)$, we define its Harish-Chandra/Selberg transform $\rho(t,n)$ by
\label{definition:harish-chandra_selberg_transform}
\[
\rho_k(t,n) = \int_{h\in SL_2(\mathbb{R})}k(I,h)y^{1/2 + it}e^{in\theta}\dfrac{dxdyd\theta}{y^2},
\]
where $y$ and $\theta$ are the Iwasawa $y$ and $\theta$-coordinates of $h$, respectively.

We are now ready to describe the general idea for bounding $K_q(g)$.

\subsection{Identifying $K_q(g)$ as the diagonal evaluation of an automorphic kernel}
\label{section:terms_on_rho}
If we manage to construct a point-pair invariant $k(g,h) = \sum_nk_n(g,h)$ with $k_n = 0$ for $n$ odd, satisfying the absolute convergence conditions discussed earlier, then the automorphic kernel $K_q(g,h)$ would have a spectral decomposition which resembles that of our $K_q(g)$ when evaluated on the diagonal. We would then identify $\rho(t,n)$ as the Harish-Chandra/Selberg transform of $k(g,h)$.

If we can prove that our chosen $k(g,h)$ attains only non-negative values, then we may bound
\[
K_q(g) = K_q(g,g) = \sum_{\gamma\in\Gamma_0(q)}k(g,\gamma g) \le \sum_{\gamma\in SL_2(\mathbb{Z})}k(g,\gamma g) = K_1(g,g) = K_1(g).
\]
The advantage of the bound $K_q(g) \le K_1(g)$ is that it is independent of $q$.

That $k(g,h)$ is real valued would follow, for example, if for all $n\in 2\mathbb{Z}$ and $g,h\in SL_2(\mathbb{R})$ one has
\[
k_n(g,h) = \overline{k_{-n}(g,h)}.
\]
The positivity would follow if we could also show that for all $g,h\in SL_2(\mathbb{R})$, $k_0(g,h) \ge 0$, and
\[
k_0(g,h) \ge \sum_{n\neq 0}|k_n(g,h)|.
\]
Replacing $g$ by $\tau g$, with $\tau\in SL_2(\mathbb{Z})$, we observe that
\[
K_q(\tau g) = \sum_{\gamma\in \Gamma_0(q)}k(\tau g,\gamma \tau g),
\]
however, since by Lemma~\ref{lemma:invariancy_lemma}, $k_n(\tau g,\gamma \tau g) = k_n(g,\tau^{-1}\gamma \tau g)$, we find that
\[
K_q(\tau g) = \sum_{\gamma\in\tau^{-1}\Gamma_0(q)\tau}k(g,\gamma g).
\]
Since $\tau^{-1}\Gamma_0(q)\tau$ is subgroup of $SL_2(\mathbb{Z})$, and since $k(g,\gamma g) \ge 0$ for all $\gamma\in SL_2(\mathbb{Z})$, we may extend our summation from $\tau^{-1}\Gamma_0(q)\tau$ to $SL_2(\mathbb{Z})$ and obtain
\[
K_q(\tau g) = \sum_{\gamma\in\tau^{-1}\Gamma_0(q)\tau}k(g,\gamma g) \le \sum_{\gamma\in\Gamma_0(1)}k(g,\gamma g) = K_1(g),
\]
which yields a bound that is not only independent on $q$, but also independent on the point $\tau g$ in the (left) orbit of $SL_2(\mathbb{Z})$ on $g$.

As an interim summary, to rectify our strategy for bounding $K_q(g)$ by $K_1(g)$, we require a test function $\rho(t,n)$ satisfying the properties:
\begin{itemize}
    \item (1) $\rho(t,n) > 0$ for all $(t,n)$ pairs of type and weight coming from eigenforms of the Casimir on $\Gamma_0(N)\char`\\SL_2(\mathbb{R})$ with $N\in\mathbb{N}$.
    \item (2) $\rho(t,n)$ is the Harish-Chandra/Selberg transform of a point-pair invariant $k(g,h)$ of the form $\sum_{n\in 2\mathbb{Z}}k_n(g,h)$, satisfying the absolute convergence property.
    \item (3) (Real valued): $k(g,h)$ attains real values.
    \item (4) (Positivity): For all $g,h\in SL_2(\mathbb{R})$ one has
    \[
    k_0(g,h) \ge \sum_{n\neq 0}|k_n(g,h)|.
    \]
    \item (5) $\rho(t,n)^{-1}$ is bounded by a polynomial in $|t|,n$.
\end{itemize}
\begin{claim}
\label{claim:rho_properties}
There exists a choice of $\rho(t,n):\mathbb{C}\times\mathbb{Z}\longrightarrow\mathbb{C}$ satisfying properties $(1)-(4)$ above, and property $(5)$ with
\[
\rho(t,n) \ll \dfrac{1}{(|t|^{2} + n^{2})^6 + 1},\quad \rho(t,n)^{-1}\ll |t|^{12} + n^{12} + 1
\]
for all pairs $(t,n)$ of type/weight coming from an eigenform of the Laplacian on $\Gamma_0(q)\char`\\SL_2(\mathbb{R})$.
\end{claim}

\textbf{Remark}: an example of such $\rho(t,n)$ is:
\[
\rho(t,n) = 
\begin{cases}
\rho_{X_0}(t)& n=0,\\
\dfrac{C}{1200}\dfrac{1}{n^{10} + 1}\rho_{X_n}(t)& n\in 2\mathbb{Z}\setminus\{0\},\\
0& \text{otherwise}.
\end{cases}
\]
Where
\[
\rho_X(t) = \dfrac{1}{t^2 + (aX)^2} + \dfrac{1}{t^2 + (bX)^2} + \dfrac{1}{t^2 + (cX)^2} + \dfrac{1}{t^2 + (dX)^2} + \dfrac{1}{t^2 + (2.5X)^2}
\]
\[
- \dfrac{1}{t^2 + X^2} - \dfrac{1}{t^2 + 4X^2} - \dfrac{1}{t^2 + 9X^2} - \dfrac{1}{t^2 + 16X^2} - \dfrac{1}{t^2 + 25X^2}.
\]
$X_n = (|n| + 2)X_0$, $X_0 = 100$, $C > 0$ is some effectively computable constant, and $a,b,c,d > 0$ are the roots of a certain quartic. For further discussion, see Appendix~\ref{appendix:constructing_rho}.

\begin{corollary}
\label{corollary:k_q_to_k_1}
For all $\tau\in SL_2(\mathbb{Z})$, one has
\[
K_q(\tau) \ll K_1(I).
\]
\end{corollary}

\subsection{Bounding $K_1(I)$}
By Corollary~\ref{corollary:k_q_to_k_1}, the proof of Theorem~\ref{theorem:bound_on_K} is reduced to:
\begin{claim}
One has
\[
K_1(I) \ll 1.
\]
\end{claim}

\begin{proof}
Recall that
\[
K_1(I) =
\sum_j\sum_n \rho(t_j,n)|v_{\lambda_j,n}(I)|^2 + \dfrac{1}{4\pi}\int_{-\infty}^{\infty}\sum_n \rho(t,n)|E_{\infty,n}(I,1/2 + it)|^2dt.
\]
The index $j$ now iterates through irreducible cuspidal representations of the Laplacian on $X_0(1)$. As for the continuous spectrum, note that we omit the outer summation over the cusps, ``$\sum_{\mathfrak{a}}$", as $X_0(1)$ has only one cusp, $\infty$.

Next, recall that (see Claim~\ref{claim:rho_properties}):
\[
\rho(t,n) \ll \dfrac{1}{(|t|^2 + n^2)^6 + 1}.
\]
The idea would be to bound each cusp form and Eisenstein series individually, and uniformly on the parameters $n$ and $t$. For cusp forms, we use a Sobolev bound due to Bernstein \& Reznikov.

\begin{claim}[cf. \cite{BERREZ}, Corollary 2.4]
Let $v_{\lambda_j,n}(I)$ be a pure weight $n$ Laplace eigenform on $X_0(1)$ of eigenvalue $\lambda_j = 1/4 + t_j^2$ and norm 1, then
\[
|v_{\lambda_j,n}(I)|^2 \ll 1 + |t_j|^2 + n^2.
\]
\end{claim}

As for the Eisenstein series, we have
\begin{claim}
Let $E_{\infty,n}(I,1/2 + it)$ be a pure weight $n$ Eisenstein series of $X_0(1)$, normalized by the constant term, then
\[
|E_{\infty,n}(I,1/2 + it)|^2 \ll \lambda(n^2\log(1 + |n|) + \lambda\log(1 + |\lambda|)).
\]
\end{claim}
For a proof of this claim, see~\cite{MUSICZEHAVI_EISENSTEIN}.

For simplicity, we bound the pointwise evaluations of both the Eisenstein series and cusp forms of type $t$ and weight $n$ at the point $I$ by a constant multiple of $1 + n^6 + |t|^6$, uniformly. Plugging these bounds back to $K_1(I)$, we obtain
\[
K_1(I) \ll \sum_j\sum_n \dfrac{1 + n^6 + |t_j|^6}{t_j^{12} + n^{12} + 1} + \int_{-\infty}^{\infty}\sum_n \dfrac{1 + n^6 + t^6}{t^{12} + n^{12} + 1}dt.
\]
Since
\[
\sum_n \dfrac{1 + n^6 + |t|^6}{t^{12} + n^{12} + 1} \ll \dfrac{1}{1 + t^4},
\]
we obtain
\[
K_1(I) \ll \sum_j\dfrac{1}{t_j^{4} + 1} + \int_{-\infty}^{\infty}\dfrac{dt}{t^{4} + 1} \ll 1 + \sum_j\dfrac{1}{t_j^{4} + 1} \ll 1 + \sum_j\dfrac{1}{\lambda_j^2 + 1}.
\]
Next, by Weyl's law for $X_0(1)$, we have $|t_j|\gg j$, so that $\sum_j\dfrac{1}{\lambda_j^2 + 1} < \infty$, so that
\[
K_1(I) \ll 1,
\]
which completes the proof.
\end{proof}
Combining the proof of the previous claim together with Corollary~\ref{corollary:k_q_to_k_1} completes the proof of Theorem~\ref{theorem:bound_on_K}.

\section{Bounding the sum $R(P)$}
\label{section:bounding_r(p)}
Recall the definition of $R(P)$.
\[
R(P) = \sum_j\sum_n\rho(t_j,n)^{-1}|<P,u_{\lambda_j,n}>|^2 + \dfrac{1}{4\pi}\sum_{\mathfrak{a}}\int_{-\infty}^{\infty}\sum_n\rho(t,n)^{-1}|<P,E_{\mathfrak{a},0}(*,1/2 + it)d>|^2dt.
\]
In this chapter, we prove:
\begin{theorem}
\label{theorem:bound_on_R}
\[
R(P) \ll_{\delta} hZ^{-19}Y^{-1-2\delta}\left(1 + h^{\frac{1}{2}}|\log Y|^2\left(Y + Y^{-1}\right)^{\frac{1}{2}}q^{-1}(h,q)^{\frac{1}{2}}\tau(hq)\right).
\]
\end{theorem}
Towards that, we first obtain an upper bound on the absolute values of the inner products. These bounds are given in terms of Fourier coefficients. We then apply the Kuznetsov/Petersson trace formulas.

\subsection{Preliminary analysis of the inner products}
\label{section:preliminaries}
Let $g\in SL_2(\mathbb{R})$, we have
\[
P(g) = \sum_{\substack{\sigma\in\Gamma_\infty \char`\\ \Gamma_0(q)}}e(hx(\sigma g))F(4\pi hy(\sigma g))G(\theta(\sigma g)),
\]
where $x(\sigma g), y(\sigma g)$ and $\theta(\sigma g)$ are the Iwasawa $x,y$ and $\theta$ coordinates of $\sigma g$, respectively.

Let $\phi_n^t$ be a pure-weight $n$ (see Definition~\ref{definition:pure_weight}) Casimir eigenfunction on $X_0(q) = \Gamma_0(q)\char`\\ SL_2(\mathbb{R})$. Unfolding,
\[
<P,\phi_n^t> = \int_{0}^{1}\int_{0}^{2\pi}\int_{0}^{\infty}e(hx)F(4\pi hy)G(\theta)\overline{\phi_n^t(x,y,\theta)}\dfrac{dxdyd\theta}{y^2}.
\]
The Fourier-Whittaker expansion of $\phi_n^t$ is given by
\[
\phi_n^t(x,y,\theta) = \sum_{m}\dfrac{\lambda_n(m)}{2\sqrt{|m|}}e(mx)e^{in\theta}W_{sgn(m)n/2,it}(4\pi |m|y).
\]
Here $sgn(\cdot)$ is the sign function, and $W_{sgn(m)n/2,it}(\cdot)$ is the Whittaker-$W$ function. Plugging this back, we obtain
\[
<P,\phi_n^t> = \dfrac{g_n\overline{\lambda_n(h)}}{2\sqrt{h}}\int_0^{\infty}y^{-2}F(4\pi hy)\overline{W_{n/2,it}(4\pi hy)}dy,
\]
where $g_n$ is the $n$'th Fourier coefficient of $G(\theta)$.

Changing variables, this becomes
\[
<P,\phi_n^t> = 2\pi\sqrt{h}g_n\overline{\lambda_n(h)}\int_0^{\infty}y^{-2}F(y)\overline{W_{n/2,it}(y)}dy.
\]
The following Lemma relates the Fourier-Whittaker coefficient $\lambda_n(h)$ to that of a minimum weight eigenform in the corresponding irreducible representation:
\begin{lemma}
\label{lemma:norm}
Let $h\in\mathbb{Z}_+$, $\phi_n^t$ and $\phi_m^t$ be a pair of pure (even) weight Casimir eigenfunctions from the same irreducible representation with $0\le |n|\le |m|$ and $nm\ge 0$. Then 
\[
|\lambda_m(h)| = |\lambda_n(h)|\cdot
\begin{cases}
\sqrt{\left|\dfrac{\Gamma(1/2 + it + |n|/2)\Gamma(1/2 - it + |n|/2)}{\Gamma(1/2 + it + |m|/2)\Gamma(1/2 - it + |m|/2)} \right|},& m\ge 0,\\\\
\sqrt{\left|\dfrac{\Gamma(1/2 + it + |m|/2)\Gamma(1/2 - it + |m|/2)}{\Gamma(1/2 + it + |n|/2)\Gamma(1/2 - it + |n|/2)} \right|},& m\le 0.
\end{cases}
\]
\end{lemma}
\begin{proof}
This is a consequence of (12) and (22) in~\cite{HARCOS-BLOMER}, alternatively, see (4.24-4.30) in~\cite{ARTIN-L-FUNCTION}. The authors would like to thank Valentin Blomer for providing them with these references.
\end{proof}

Our analysis is split between the cases:
\begin{itemize}
    \item Even Principal Series of small type parameter $(|t|\le 1)$.
    \item Discrete Series $(it\in \mathbb{Z}/2)$.
    \item Even Principal Series of large type parameter $(t\ge 1)$.
\end{itemize}

\textbf{Remark}: in the case of small $t$ parameter, we include both the case where $0\le t\le 1$ and the case where $t$ is exceptional, i.e. $0\le it\le 1/4$. The positivity of the Casimir, together with $\lambda = 1/4 + t^2$, imply that the exceptional eigenvalues satisfy the inequality $0\le it\le 1/2$. The refinement $0\le it\le 1/4$ comes from Selberg's $3/16$ Theorem.

\subsection{The Even principal/Eisenstein series for $|t|\le 1$}
\label{section:small_parameter}
\begin{theorem}
\label{theorem:small_parameter}
Let $\phi_n^t$ be a pure weight $n$ Casimir eigenfunction of type $|t| \le 1$, and let $0 < \delta < 1/4$. We have
\[
|<P,\phi_n^t>| \ll_{\delta} \sqrt{h}|g_n||\lambda_0(h)|(1 + |n/2|)Y^{-1/2 - \delta}(Y^{it} + Y^{-it} + 3).
\]
\end{theorem}
We showed that:
\[
|<P,\phi_n^t>| = 2\pi\sqrt{h}|g_n||\lambda_n(h)|\left|\int_0^{\infty}y^{-2}F(y)\overline{W_{n/2,it}(y)}dy\right|.
\]
Lemma~\ref{lemma:norm} gives:
\[
|<P,\phi_n^t>| = 
2\pi\sqrt{h}|g_n||\lambda_0(h)|\cdot \sqrt{\left|\dfrac{\Gamma(1/2 + it)\Gamma(1/2 - it)}{\Gamma(1/2 + it + |n/2|)\Gamma(1/2 - it + |n/2|)}\right|}^{sgn(n)}
\left|\int_0^{\infty}y^{-2}F(y)\overline{W_{n/2,it}(y)}dy\right|.
\]
Theorem~\ref{theorem:small_parameter} reduces to
\begin{claim}
\[
\sqrt{\left|\dfrac{\Gamma(1/2 + it)\Gamma(1/2 - it)}{\Gamma(1/2 + it + |n/2|)\Gamma(1/2 - it + |n/2|)}\right|}^{sgn(n)}\left|\int_0^{\infty}y^{-2}F(y)\overline{W_{n/2,it}(y)}dy\right| \ll 
(1 + |n/2|)Y^{-1/2 - \delta}(Y^{it} + Y^{-it} + 3).
\]
\end{claim}
To complete the analysis, we apply
\begin{claim}
\[
\sqrt{\left|\dfrac{\Gamma(1/2 + it)\Gamma(1/2 - it)}{\Gamma(1/2 + it + |n/2|)\Gamma(1/2 - it + |n/2|)}\right|}^{sgn(n)}|W_{n/2,it}(y)| \ll_{\delta} (|n/2| + |it| + 1)y^{1/2 - |\Re(it)| - \delta}.
\]
\end{claim}
\begin{proof}
See Bruggeman and Motohashi~\cite[4.3]{MOTOHASHI}. Their bound is valid for $y \longrightarrow 0^+$, which is our case.
\end{proof}
Applying the triangle inequality and the Bruggeman-Motohashi bound completes the proof of Theorem~\ref{theorem:small_parameter}.

\subsection{The discrete series}
As a first step, we explain why we may restrict our analysis to the holomorphic series.
\subsubsection{Reducing to the holomorphic series case}
Let $V$ denote the irreducible representation containing $\phi_n^t$. $V$ has a minimal $K$-type, which is a pair of pure-weight one-dimensional Casimir eigenspaces $V_k$ and $V_{-k}$, such that $V_{l} = 0$ for all $|l| < |k|$; $k\neq 0$.

Let us denote by $\phi_k^t$ and $\phi_{-k}^t$ a pair of unit vectors of $V_k$ and $V_{-k}$, respectively. Let $n\in 2\mathbb{Z}_+$. It is well known that $\phi_{n+k}^t = \overline{\phi_{-n-k}^t}$, which implies that $|<P,\phi_{-n-k}^t>| = |<P,\phi_{n+k}^t>|$. Therefore, the contribution of the holomorphic series and the anti-holomorphic series to $R(P)$ is equal.

\subsubsection{The holomorphic series}
\begin{theorem}[Holomorphic series]
\label{theorem:discrete_series}
With all notations as above, $0\le m\in2\mathbb{Z}$, we have
\[
|<P,\phi_{k+m}^t>| \ll \sqrt{\Gamma(k)}\sqrt{h}|g_{k+m}||\lambda_k(h)|Y^{-1/2}(k + m)^2.
\]
\end{theorem}
The proof of this theorem is more intricate than Theorem~\ref{theorem:small_parameter}, as the Bruggeman-Motohashi isn't sharp enough for our application.

To prove this theorem one needs to appeal to `The Mellin-Barnes Integral" representation of the Whittaker-W function, see~\cite[7.621.11]{RYZHIK}, which is essentially the inverse Mellin transform of its Mellin transform.
\begin{claim}
\label{claim:mellin_barnes}
Let $s=\sigma + ir$, then for $\sigma$ sufficiently large (such that all poles of the integrand are to the left of the line $\Re(s) = \sigma$), one has
\[
W_{m,it}(y)e^{-y/2} = \int_{\Re(s)=\sigma}\dfrac{\Gamma(1/2+s-it)\Gamma(1/2+s+it)}{\Gamma(1+s-m)}y^{-s}\dfrac{ds}{2\pi i}.
\]
\end{claim}
For a proof of Theorem~\ref{theorem:discrete_series}, see~\cite{SAARTHESIS}.

\subsection{Analyzing the even principal/Eisenstein series for $t\ge 1$}
\label{section:large_parameter}
\begin{theorem}[Eisenstein series of real parameter $t \ge 1$]
\label{theorem:large_parameter}
Let $\phi_n^t$ be a pure-weight $n \in 2\mathbb{Z}$ Casimir eigenfunction with $|t|\ge 1$, one has
\[
|<P,\phi_n^t>| \ll \sqrt{h}|g_n||\lambda_0(h)|\dfrac{Y^{-1/2}}{\sqrt{\cosh(\pi t)}(1 + t^{12})}I_0(2\sqrt{(|n/2| + 1)Y}),
\]
where
\[
I_{\alpha}(x) := \sum_{k=0}^{\infty}\dfrac{(x/2)^{2k+\alpha}}{k!\Gamma(k + \alpha + 1)}
\]
is the modified Bessel function of the first kind,
\end{theorem}
Similar to the case of the discrete series, the Bruggeman-Motohashi bound isn't strong enough for our application. Instead, we use the Whittaker-W series expansion. 

\begin{proof}
In Section~\ref{section:preliminaries}, it has been established that 
\[
|<P,\phi_n^t>| = 2\pi\sqrt{h}|g_n||\lambda_n(h)|\left|\int_0^{\infty}y^{-2}F(y)\overline{W_{n/2,it}(y)}dy\right|.
\]
An immediate implication of Lemma~\ref{lemma:norm} yields
\[
|<P,\phi_n^t>| =
2\pi\sqrt{h}|g_n||\lambda_0(h)|\cdot \sqrt{\left|\dfrac{\Gamma(1/2 + it)\Gamma(1/2 - it)}{\Gamma(1/2 + it + |n/2|)\Gamma(1/2 - it + |n/2|)}\right|}^{sgn(n)}\left|\int_0^{\infty}y^{-2}F(y)\overline{W_{n/2,it}(y)}dy\right|.
\]
Theorem~\ref{theorem:large_parameter} reduces to:
\begin{claim}
\label{claim:large_parameter_integral}
With $t,n$ as in Theorem~\ref{theorem:large_parameter}, one has
\[
\left|\dfrac{\Gamma(1/2 + it)}{\Gamma(1/2 + it + |n/2|)}\right|^{sgn(n)}\left|\int_0^{\infty}y^{-2}F(y)\overline{W_{n/2,it}(y)}dy\right| 
\ll \dfrac{Y^{-1/2}}{\sqrt{\cosh(\pi t)}(1 + t^{12})}\sum_{k=0}^{\infty}\dfrac{(|n/2| + 1)^kY^k}{k!^2}.
\]
\end{claim}

\begin{proof}
\begin{claim}[See~\cite{WOLFRAM}]
\label{claim:power_series}
\[
W_{m,it}(y) =
e^{-y/2}\bigg(\sum_{k=0}^{\infty}
\dfrac{1}{k!}
\dfrac{\Gamma(2it)}{\Gamma(1/2 + it - m)}
\dfrac{\Gamma(1/2 - it - m + k)}{\Gamma(1/2 - it - m)}
\dfrac{\Gamma(1 - 2it)}{\Gamma(1 - 2it + k)}
y^{1/2+k-it}
\]
\[
+ 
\dfrac{1}{k!}
\dfrac{\Gamma(-2it)}{\Gamma(1/2 - it - m)}
\dfrac{\Gamma(1/2 + it - m + k)}{\Gamma(1/2 + it - m)}
\dfrac{\Gamma(1 + 2it)}{\Gamma(1 + 2it + k)}
y^{1/2+k+it}\bigg).
\]
\end{claim}
Substituting into $\left|\dfrac{\Gamma(1/2 + it)}{\Gamma(1/2 + it + |n/2|)}\right|^{sgn(n)}\left|\int_0^{\infty}y^{-2}F(y)\overline{W_{n/2,it}(y)}dy\right|$, and applying the triangle inequality, we obtain the upper bound:
\[
\left|\dfrac{\Gamma(1/2 + it)}{\Gamma(1/2 + it + |n/2|)}\right|^{sgn(n)}\times
\]
\[
\Bigg(\sum_{k=0}^{\infty}\left|
\dfrac{1}{k!}
\dfrac{\Gamma(2it)}{\Gamma(1/2 + it - n/2)}
\dfrac{\Gamma(1/2 - it - n/2 + k)}{\Gamma(1/2 - it - n/2)}
\dfrac{\Gamma(1 - 2it)}{\Gamma(1 - 2it + k)}
\int_{y=0}^{\infty}e^{-y/2}F(y)y^{1/2+k-it}\dfrac{dy}{y^2}\right|
\]
\[
+ 
\left|\dfrac{1}{k!}
\dfrac{\Gamma(-2it)}{\Gamma(1/2 - it - n/2)}
\dfrac{\Gamma(1/2 + it - n/2 + k)}{\Gamma(1/2 + it - n/2)}
\dfrac{\Gamma(1 + 2it)}{\Gamma(1 + 2it + k)}
\int_{y=0}^{\infty}e^{-y/2}F(y)y^{1/2+k+it}\dfrac{dy}{y^2}\right|\Bigg).
\]
The two summands inside the parenthesis are conjugate, hence the left hand side of Claim~\ref{claim:large_parameter_integral} is bounded above by a constant multiple of
\[
\Bigg|\left(\dfrac{\Gamma(1/2 + it)}{\Gamma(1/2 + it + |n/2|)}\right)^{sgn(n)}
\sum_{k=0}^{\infty}\dfrac{1}{k!}
\dfrac{\Gamma(2it)}{\Gamma(1/2 + it - n/2)}
\dfrac{\Gamma(1/2 - it - n/2 + k)}{\Gamma(1/2 - it - n/2)}
\dfrac{\Gamma(1 - 2it)}{\Gamma(1 - 2it + k)}
\]
\[
\times
\int_{y=0}^{\infty}e^{-y/2}F(y)y^{1/2+k-it}\dfrac{dy}{y^2}\Bigg|.
\]

\begin{claim}
\label{claim:large_parameter_inner_integral}
One has,
\[
\left|\int_{y=0}^{\infty}e^{-y/2}F(y)y^{1/2+k-it}\dfrac{dy}{y^2}\right| \ll 
\dfrac{Y^{-1/2 + k}}{(1 + t^{12})}.
\]
\end{claim}
\begin{proof}
This is a consequence of integration by parts 12 times,
using the fact that $F(y)$ is supported on $[Y,2Y]$, where $0<Y\ll 1$, and $|F^{(j)}|\ll_j Y^{-j}$. See Definition~\ref{definition:function_F}.
\end{proof}
Claim~\ref{claim:large_parameter_integral} reduces to:
\begin{claim}
\label{claim:large_parameter_gamma_sum}
\small{
\[
\left|\left(\dfrac{\Gamma(1/2 + it)}{\Gamma(1/2 + it + |n/2|)}\right)^{sgn(n)}
\sum_{k=0}^{\infty}\dfrac{Y^k}{k!}
\dfrac{\Gamma(2it)}{\Gamma(1/2 + it - n/2)}
\dfrac{\Gamma(1/2 - it - n/2 + k)}{\Gamma(1/2 - it - n/2)}
\dfrac{\Gamma(1 - 2it)}{\Gamma(1 - 2it + k)}\right|
\ll \dfrac{1}{\sqrt{\cosh(\pi t)}}\sum_{k=0}^{\infty}\dfrac{(|n/2| + 1)^kY^k}{k!^2}.
\]
}
\end{claim}
\begin{proof}
We shift our focus towards the summand,
\[
\left|\dfrac{\Gamma(1/2 + it)}{\Gamma(1/2 + it + |n/2|)}\right|^{sgn(n)}\left|
\dfrac{\Gamma(2it)}{\Gamma(1/2 + it - n)}
\dfrac{\Gamma(1/2 - it - n/2 + k)}{\Gamma(1/2 - it - n/2)}
\dfrac{\Gamma(1 - 2it)}{\Gamma(1 - 2it + k)}
\right|.
\]
We prove:
\begin{claim}
\label{claim:large_parameter_easy_gamma_ratio}
\[
\left|\dfrac{\Gamma(1/2 + it)}{\Gamma(1/2 + it + |n/2|)}\right|^{sgn(n)}\left|
\dfrac{\Gamma(2it)}{\Gamma(1/2 + it - n/2)}\right| \ll \dfrac{1}{\sqrt{\cosh(\pi t)}}.
\]
\end{claim}
\begin{proof}
We split the analysis into two cases with respect to whether $n$ is positive or negative.

\textbf{For $\mathbf{n\ge0}$}: the reflection formula shows that
\[
\left|\dfrac{1}{\Gamma(1/2 + it + n/2)\Gamma(1/2 + it - n/2)}\right| \ll \cosh(\pi t),
\]
whereas, combining the identities
\[
|\Gamma(2it)|^2 = \dfrac{\pi}{2t\sinh(2\pi t)},\quad |\Gamma(1/2 + it)|^2 = \dfrac{\pi}{\cosh(\pi t)},
\]
see~\cite[8.332.1-2]{RYZHIK}, and multiplying these identities with the above upper bound completes the proof.

\textbf{For $\mathbf{n<0}$}: the factors $\Gamma(1/2 + it - n/2)$ and $\Gamma(1/2 + it + |n/2|)$ cancel each other out and we are left with bounding $\left|\dfrac{\Gamma(2it)}{\Gamma(1/2 + it)}\right|$, which follows immediately from the above identities.
\end{proof}

Claim~\ref{claim:large_parameter_gamma_sum} reduces to:
\begin{claim}
\label{claim:large_param_gamma_sum_equivalent_form}
\[
\sum_{k=0}^{\infty}
\left|\dfrac{\Gamma(1/2 - it - n/2 + k)}{\Gamma(1/2 - it - n/2)}
\dfrac{\Gamma(1 - 2it)}{\Gamma(1 - 2it + k)}\right|\cdot \dfrac{Y^k}{k!} \ll \sum_{k=0}^{\infty}\dfrac{(|n/2| + 1)^kY^k}{k!^2}.
\]
\end{claim}
\begin{proof}
Denote the gamma factor by $C(k)$, i.e.
\[
C(k) = \dfrac{\Gamma(1/2 - it - n/2 + k)\Gamma(1 - 2it)}{\Gamma(1/2 - it - n/2)\Gamma(1 - 2it + k)}.
\]
We are given the task of bounding the sum
\[
\sum_{k=0}^{\infty}|C(k)|\dfrac{Y^k}{k!}.
\]
We approach this goal in steps, starting with a naive bound on $|C(k)|$.

\begin{claim}
For all $k\in\mathbb{Z}_{+}$, one has
\[
|C(k)| \le \prod_{j=1}^k(1 + \dfrac{|n/2|+1}{j}).
\]
\end{claim}

\begin{proof}
Clearly $C(0) = 1$. The recursive structure of the gamma function yields
\[
C(k) = C(k-1)\cdot \dfrac{- it - n/2 + k - 1/2}{-2it + k},
\]
which gives
\[
C(k) = \prod_{j=1}^k\dfrac{- it - n/2 + j - 1/2}{- 2it + j}.
\]
Clearly,
\[
|\dfrac{- it - n/2 + j - 1/2}{- 2it + j}|\le 1 + \dfrac{|n/2| + 1}{j}.
\]
Applying this bound term by term on the formula for $C(k)$ completes the proof.
\end{proof}
Denote
\[
D(k) := \prod_{j=1}^k(1 + \dfrac{|n/2|+1}{j}).
\]
Claim~\ref{claim:large_param_gamma_sum_equivalent_form} reduces to:
\begin{claim}
\label{claim:large_param_gamma_sum_equivalent_form_number_two}
\[
\sum_{k=0}^{\infty}\dfrac{D(k)Y^k}{k!} \ll \sum_{k=0}^{\infty}\dfrac{|n/2|^kY^k}{k!^2}.
\]
\end{claim}

\begin{proof}
\begin{claim}[Step 1]
For $Y < \dfrac{1}{2e}$ one has
\[
\sum_{k=0}^{\infty}\dfrac{D(k)Y^k}{k!} \ll \sum_{k=0}^{|n/2|}\dfrac{D(k)Y^k}{k!}.
\]
\end{claim}
\begin{proof}
For all $k\in\mathbb{Z}_+$, $D(k) \le (|n/2| + 2)^k$, and in particular $D(|n/2|) \le (|n/2| + 2)^{|n/2|}$. Also, for $j > |n/2|$, we have $1 + \frac{|n/2| + 1}{j} \le 2$, and therefore, for all $k > |n/2|$ we have
\[
\dfrac{D(k)Y^k}{k!}\le \dfrac{D(|n/2|)Y^k\cdot\prod_{j=|n/2|+1}^k(1 + \frac{|n/2| + 1}{j})}{k!}
\le \dfrac{(|n/2|+2)^{|n/2|}Y^k2^{k-|n/2|}}{|n/2|!} \ll \left(2eY\right)^k,
\]
which implies that for $Y < \dfrac{1}{2e}$:
\[
\sum_{k=|n/2|+1}^{\infty}\dfrac{D(k)Y^k}{k!} \ll 1.
\]
Noting that $\dfrac{D(0)Y^0}{0!} = 1$, the claim follows.
\end{proof}
\begin{claim}[Step 2]
For $Y < \dfrac{1}{2e^2}$ one has
\[
\sum_{k=0}^{|n/2|}\dfrac{D(k)Y^k}{k!} \ll \sum_{k=0}^{\lceil\sqrt{|n/2|}\rceil}\dfrac{D(k)Y^k}{k!}.
\]
\end{claim}
\begin{proof}
Note that for all $0\le j\le |n/2|$, $\dfrac{|n/2|+1}{j}\ge 1$, which implies that
\[
\forall\quad 0\le k\le |n/2|:\quad D(k) = \prod_{j=1}^k(1 + \dfrac{|n/2|+1}{j})\le \prod_{j=1}^k(2\cdot\frac{|n/2|+1}{j}) = 2^k\dfrac{(|n/2|+1)^k}{k!}.
\]
Let $\lceil\sqrt{|n/2|}\rceil < k \le |n/2|$, we have
\[
\dfrac{D(k)Y^k}{k!} \le \dfrac{2^k(|n/2|+1)^kY^k}{(k!)^2} \ll \left(\dfrac{2(|n/2|+1)e^2Y}{k^2}\right)^k
\ll \left(2e^2Y\right)^k,
\]
which implies that for $Y < \dfrac{1}{2e^2}$:
\[
\sum_{k=\lceil\sqrt{|n/2|}\rceil + 1}^{|n/2|}\dfrac{D(k)Y^k}{k!} \ll 1,
\]
which completes the proof.
\end{proof}
\begin{claim}[Step 3]
\[
\sum_{k=0}^{\lceil\sqrt{|n/2|}\rceil}\dfrac{D(k)Y^k}{k!} \ll \sum_{k=0}^{\lceil\sqrt{|n/2|}\rceil}\dfrac{(|n/2| + 1)^kY^k}{k!^2}.
\]
\end{claim}
\begin{proof}
We require a finer bound on $D(k)$, valid for $0\le k \le \lceil\sqrt{|n/2|}\rceil$. The idea is to estimate $D(k) = \prod_{j=1}^k(1 + \dfrac{|n/2|+1}{j})$ with $\prod_{j=1}^k\dfrac{|n/2|+1}{j}$. 

Consider the ratio
\[
\dfrac{\prod_{j=1}^k(1 + \dfrac{|n/2|+1}{j})}{\prod_{j=1}^k\dfrac{|n/2|+1}{j}} = \prod_{j=1}^k(1 + \dfrac{j}{|n/2|+1})\le
e^{\frac{1}{|n/2| + 1}\sum_{j=1}^kj}.
\]
Since $0\le k\le \lceil\sqrt{|n/2|}\rceil$,
\[
\dfrac{1}{|n/2| + 1}\sum_{j=1}^kj \ll 1 \implies D(k) \ll \prod_{j=1}^k\dfrac{|n/2|+1}{j}.
\]
Therefore, for $0\le k \le \lceil\sqrt{|n/2|}\rceil$,
\[
D(k) \ll \dfrac{(|n/2| + 1)^k}{k!},
\]
which completes the proof.
\end{proof}
Recall that under our assumptions on $h$, see Theorem~\ref{theorem:point-wise_bound}, $h < N^{1/3}$. $Y$ is defined as $4\pi h/N$, hence, assuming $N$ is large enough, the assumptions of the above claims are satisfied, so that Claim~\ref{claim:large_param_gamma_sum_equivalent_form_number_two} follows from extending our summation indefinitely.
\end{proof}
Claim~\ref{claim:large_param_gamma_sum_equivalent_form_number_two} implies Claim~\ref{claim:large_param_gamma_sum_equivalent_form}.
\end{proof}
Claim~\ref{claim:large_parameter_gamma_sum} follows from Claim~\ref{claim:large_param_gamma_sum_equivalent_form} and Claim~\ref{claim:large_parameter_easy_gamma_ratio}.
\end{proof}
Claim~\ref{claim:large_parameter_integral} follows from Claim~\ref{claim:large_parameter_gamma_sum} and Claim~\ref{claim:large_parameter_inner_integral}.
\end{proof}
Theorem~\ref{theorem:large_parameter} follows from Claim~\ref{claim:large_parameter_integral}.
\end{proof}

\subsection{Transforming the bound on $R(P)$ - Inner summation}
\label{section:transforming_the_bound_on_R}
Recall the definition of $R$,
\[
R(P) = \sum_j\sum_n\rho(n,t)^{-1}|<P,u_n^{t_j}>|^2 + \int_{-\infty}^{\infty}\sum_n\rho(n,t)^{-1}|<P,E_n(*,1/2 + it)>|^2dt,
\]
where $\rho(n,t)^{-1}$ is our test function, see Claim~\ref{claim:rho_properties}.

Over the past section we bounded $|<P,\phi_n^t>|$. Our analysis was conducted separately for different types of representations. In this section, we combine these bounds in order to construct an upper bound on the inner summations: $\sum_n\rho(n,t)^{-1}|<P,\phi_n^t>|^2$. Similar to the previous section, our analysis is split between the cases:
\begin{itemize}
    \item Even Principal Series of small type parameter $(|t|\le 1)$.
    \item Even Principal Series of large type parameter $(t\ge 1)$.
    \item Discrete Series $(it\in \mathbb{Z}/2)$
\end{itemize}

\subsubsection{Summing the spectral estimates for $|t| \le 1$}
\begin{claim}
\label{claim:reformulation_small_parameter}
Let $\{\phi_n^t\}_{n\in 2\mathbb{Z}}$ be an even principal series of type $t$, $|t| \le 1$, we have:
\[
\sum_n\rho(n,t)^{-1}|<P,\phi_n^t>|^2 \ll h|\lambda_0(h)|^2Z^{-13}Y^{-1-2\delta}(Y^{2it} + Y^{-2it} + 3).
\]
\end{claim}

\begin{proof}
Plugging the bound from Theorem~\ref{theorem:small_parameter}, while recalling that because $|t| \le 1$, we have $\rho^{-1}(n,t) \ll 1 + |n|^{12}$, we obtain
\[
\sum_{n\in2\mathbb{Z}}\rho(n,t)^{-1}|<P,\phi_n^t>|^2 \ll h|\lambda_0(h)|^2Y^{-1-2\delta}(Y^{2it} + Y^{-2it} + 3)\left(\sum_{n\in\mathbb{Z}}|g_n|^2(1 + |n|^{14}) \right).
\]
As for $\sum_{n\in\mathbb{Z}}|g_n|^2(1 + |n|^{14})$: 
\[
\sum_{n\in\mathbb{Z}}|g_n|^2 = ||G||^2\text{ and } \sum_{n\in\mathbb{Z}}|g_n|^2n^{14} = ||G^{(7)}||^2.
\]
Since $||G^{(7)}||^2 \ll Z^{-13}$, $||G|| \ll 1$, and $Z \ll 1$, see Definition~\ref{definition:function_G}, 
we find that
\[
\sum_{n\in\mathbb{Z}}|g_n|^2(1 + |n|^{14}) \ll Z^{-13},
\]
which completes the proof.
\end{proof}

\subsubsection{Summing the spectral estimates for $t \ge 1$}
\begin{claim}
\label{claim:reformulation_large_parameter}
Let $\{\phi_n^t\}_{n\in 2\mathbb{Z}}$ be an even principal series of type $t$, $t \ge 1$, then:
\[
\sum_n\rho(n,t)^{-1}|<P,\phi_n^t>|^2 \ll \dfrac{h|\lambda_0(h)|^2Z^{-11}Y^{-1}}{\cosh(\pi t)(1 + t^{12})}.
\]
\end{claim}
The bound obtained in Theorem~\ref{theorem:large_parameter} holds in our case. Recalling that $\rho^{-1}(n,t) \ll 1 + t^{12} + |n|^{12} \le (1 + t^{12})(1 + |n|^{12})$, we have
\[
\sum_n\rho(n,t)^{-1}|<P,\phi_n^t>|^2 \ll \dfrac{h|\lambda_0(h)|^2Y^{-1}}{\cosh(\pi t)(1 + t^{12})}
\left(\sum_{n\in 2\mathbb{Z}}|g_n|^2(1 + |n|^{12})I_0(2\sqrt{(|n|/2 + 1)Y})^2 \right).
\]
Claim~\ref{claim:reformulation_large_parameter} reduces to:
\begin{claim}
The following estimate holds.
\[
\sum_{n\in 2\mathbb{Z}}|g_n|^2(1 + |n|^{12})I_0(2\sqrt{(|n|/2 + 1)Y})^2
 \ll Z^{-11}.
\]
\end{claim}
\begin{proof}
We split the sum into two parts.
\[
\sum_{n\in 2\mathbb{Z}}|g_n|^2(1 + |n|^{12})I_0(2\sqrt{(|n|/2 + 1)Y})^2 = \sum_{n\in 2\mathbb{Z}: |n|Y \le 1} + \sum_{n\in 2\mathbb{Z}: |n|Y > 1}.
\]
For $|n|Y \le 1$, we have $I_0(2\sqrt{(|n|/2 + 1)Y})^2 \ll 1$. Therefore,
\[
\sum_{n\in 2\mathbb{Z}: |n|Y \le 1} \ll \sum_{n\in 2\mathbb{Z}: |n|Y \le 1}|g_n|^2(1 + |n|^{12}) \le ||G||^2 + ||G^{(6)}||^2 \ll Z^{-11}.
\]
As for $\sum_{n\in 2\mathbb{Z}: |n|Y > 1}$, we bound $I_0(|x|) \le e^{|x|}$, and $|g_n| \ll e^{-2\sqrt{\pi|n|Z}}$ (see Definition~\ref{definition:function_G}). Plugging back, we obtain:
\[
\sum_{n\in 2\mathbb{Z}: |n|Y > 1}
\ll
\sum_{n\in \mathbb{N}: nY > 1}n^{12}e^{4\sqrt{\pi n}(\sqrt{Y} - \sqrt{Z})}.
\]
The function $x^{12}e^{4\sqrt{\pi x}(\sqrt{Y} - \sqrt{Z})}$ is positive and monotonically decreasing for $xY > 1$. Therefore,
\[
\sum_{n\in \mathbb{N}: nY > 1} \ll (1/Y + 1)^{12}e^{4\sqrt{\pi (1/Y + 1)}(\sqrt{Y} - \sqrt{Z})} + \int_{x = 1/Y}^{\infty}x^{12}e^{4\sqrt{\pi x}(\sqrt{Y} - \sqrt{Z})}dx.
\]
Because $Y < N^{-2/3}$, while $Z > N^{-1/3}$, we find that both terms vanish as $N$ tends to infinity, and therefore $\sum_{n\in 2\mathbb{Z}: |n|Y > 1} \ll 1$. Combining the two bounds, we find that
\[
\sum_{n\in 2\mathbb{Z}}|g_n|^2(1 + |n|^{12})I_0(2\sqrt{(|n|/2 + 1)Y})^2
 \ll Z^{-11},
\]
which completes the proof.
\end{proof}
\subsubsection{Summing the spectral estimates for the discrete series}
\begin{claim}
\label{claim:reformulation_discrete_series}
Let $\{\phi_{\pm (k+n)}^t\}_{n\in 2\mathbb{Z}_{\ge 0}}$ be a discrete series with min-weight base vector $\phi_{\pm k}^t$, $k\in\mathbb{Z}_{+}$, then
\[
\sum_{n\ge 0}\rho(\pm(k+n),t)^{-1}|<P,\phi_{\pm(k+n)}^t>|^2 \ll \dfrac{h\Gamma(k)|\lambda_k(h)|^2Z^{-19}Y^{-1}}{k^4}.
\]
\end{claim}
\begin{proof}
Plugging the bound from Theorem~\ref{theorem:discrete_series}, while recalling that $it = \dfrac{k-1}{2}$, and hence $\rho^{-1}(\pm(k+n),t) \ll (k+n)^{12}$, we obtain:
\[
\sum_{n\ge 0}\rho(\pm(k+n),t)^{-1}|<P,\phi_{\pm(k+n)}^t>|^2 \ll \dfrac{h\Gamma(k)|\lambda_k(h)|^2Y^{-1}}{k^4} \sum_{n\ge 0}|g_{k+n}|^2(k + n)^{20}.
\]
Since $\sum_{n\ge 0}|g_{k+n}|^2(k + n)^{20} \le ||G^{(10)}||^2 \ll Z^{-19}$, the claim follows.
\end{proof}

\subsection{Transforming the bound on $R(P)$ - Outer summation}
\label{section:transforming_the_bound_on_R_outer_summation}
Denote by $R_E(P)$ and $R_D(P)$ the contributions of the even principal series and the discrete series to $R(P)$, respectively. In this section, we prove:
\begin{claim}
\label{claim:r_k_bound}
\[
R_E(P) \ll hZ^{-13}Y^{-1-2\delta}\left(1 + \log^2 (h/Y)h^{\frac{1}{2}}\left(Y + Y^{-1}\right)^{\frac{1}{2}}q^{-1}(h,q)^{\frac{1}{2}}\tau(hq)\right).
\]
\end{claim}
And
\begin{claim}
\label{claim:r_p_bound}
\[
R_D(P) \ll hZ^{-19}Y^{-1}\left(1 + h^{1/2}\log^2 (h/Y)q^{-1}(h,q)^{\frac{1}{2}}\tau(hq)\right).
\]
\end{claim}

\subsubsection{Applying the Kuznetsov trace formula}
We define the even $C^{\infty}$ function
\[
g_K(t) = \dfrac{Y^{2it} + Y^{-2it} + 3}{\cosh(\pi t)(1 + t^{12})}.
\]
Combining the bounds obtained in Claim~\ref{claim:reformulation_small_parameter} and Claim~\ref{claim:reformulation_large_parameter}, we obtain:
\label{claim:even_series_bound}
\[
R_E(P) \ll_{\delta} hZ^{-13}Y^{-1-2\delta}\left(\sum_jg_K(t_j)|\lambda_{0,j}(h)|^2 + 
\dfrac{1}{4\pi}\sum_{\mathfrak{a}}\int_{-\infty}^{\infty}g_K(t)|\lambda_{0,\mathfrak{a},t}(h)|^2dt\right).
\]
Denote by $Q_E$ the expression
\[
Q_E := \sum_jg_K(t_j)|\lambda_{0,j}(h)|^2 + 
\dfrac{1}{4\pi}\sum_{\mathfrak{a}}\int_{-\infty}^{\infty}g_K(t)|\lambda_{0,\mathfrak{a},t}(h)|^2dt.
\]
Claim~\ref{claim:r_k_bound} reduces to:
\begin{claim}
\label{claim:applying_kuznetsov}
\[
Q_E \ll 1 + \log^2 (h/Y)h^{\frac{1}{2}}\left(Y + Y^{-1}\right)^{\frac{1}{2}}q^{-1}(h,q)^{\frac{1}{2}}\tau(hq).
\]
\end{claim}
\begin{proof}
The Kuznetsov trace formula for $\Gamma_{\infty}\char`\\\Gamma_0(q)$, see~\cite{DFI} for reference, implies that
\[
Q_E = \dfrac{1}{\pi^2}\int_{-\infty}^{\infty}t\sinh(\pi t)g_K(t)dt + \sum_{c = 0(q)}\dfrac{1}{c}S(h,h;c)g_K^{+}(\dfrac{4\pi h}{c}),
\]
where $S(h,h;c)$ are Kloosterman sums, $g_K^{+}(x) = \dfrac{2i}{\pi}\int_{-\infty}^{\infty}J_{2it}(x)g_K(t)tdt$,
and $J_{\nu}(y)$ is the $J$-Bessel function, see~\cite{RYZHIK}[8.40] for reference. Clearly,
\[
\int_{-\infty}^{\infty}t\sinh(\pi t)g_K(t)dt \ll 1.
\]
In order to estimate $g_K^{+}(x)$, we make the change of variables $s = 1/2 + it$ and push the contour of integration to the right, towards the line $s = \sigma + it$, with $1/2 \le \sigma < 1$. Using the Cauchy Residue Theorem, since we are not picking up any poles and the integrand vanishes at infinity, we have
\[
g_K^{+}(x) = \dfrac{2}{\pi i}\int_{s = \sigma + it}J_{2s-1}(x)g_K((s-\dfrac{1}{2})i)(s-\dfrac{1}{2})ds.
\]
We have $|J_{2s-1}(x)| \ll e^{\pi|s|}x^{2\sigma - 1}$, see~\cite{RYZHIK}[8.411.6], and $
|g_K((s-\dfrac{1}{2})i)| \ll \left(Y^{2\sigma - 1} + Y^{1 - 2\sigma} + 3\right)|s|^{-12}e^{-\pi|s|}$. Therefore
\[
g_K^{+}(x) \ll x^{2\sigma - 1}\left(Y^{2\sigma - 1} + Y^{1 - 2\sigma} + 3\right).
\]
Plugging back, we obtain
\[
Q_E \ll 1 + h^{2\sigma - 1}\left(Y^{2\sigma - 1} + Y^{1 - 2\sigma} + 3\right)\sum_{c=0(q)}c^{-2\sigma}|S(h,h;c)|.
\]
Applying Weil's bound on Kloosterman sums, one has $|S(h,h;c)| \le (h,c)^{\frac{1}{2}}c^{\frac{1}{2}}\tau(c)$. To ensure the convergence of the sum, we restrict our choice of $\sigma$ to the interval $(\dfrac{3}{4}, 1)$. We obtain (cf.~\cite[p. 431]{DFI}):
\[
\sum_{c=0(q)}c^{\frac{1}{2}-2\sigma}(h,c)^{\frac{1}{2}}\tau(c) \ll (\sigma - \dfrac{3}{4})^{-2}q^{\frac{1}{2} - 2\sigma}(h,q)^{\frac{1}{2}}\tau(hq).
\]
Choosing $\sigma = \dfrac{3}{4} + \dfrac{1}{\log h/Y}$ and plugging back, we obtain
\[
Q_E \ll 1 + \log^2 (h/Y)h^{\frac{1}{2}}\left(Y + Y^{-1}\right)^{\frac{1}{2}}q^{-1}(h,q)^{\frac{1}{2}}\tau(hq),
\]
which completes the proof.
\end{proof}

\subsubsection{Applying the Petersson trace formula}
By Claim~\ref{claim:reformulation_discrete_series}, we have
\[
R_D(P) \ll hZ^{-19}Y^{-1}\sum_{k\in 2\mathbb{N}}\dfrac{\Gamma(k)}{k^4}\sum_{f\in\mathcal{F}_k}|\lambda_f(h)|^2,
\]
where $\sum_{f\in\mathcal{F}_k}$ runs over an orthonormal basis of the complex vector spaces of min-weight $k$ vectors of irreducible (holomorphic) discrete series representations. Denote by $Q_P(k)$ the expression
\[
Q_P(k) := \Gamma(k)\sum_{f\in\mathcal{F}_k}|\lambda_f(h)|^2,
\]
such that $R_D(P) \ll hZ^{-19}Y^{-1}\sum_{k\in 2\mathbb{N}}\dfrac{1}{k^4}Q_P(k)$. Claim~\ref{claim:r_p_bound} reduces to:
\begin{claim}
\label{claim:applying_petersson}
\[
Q_P(k) \ll k\left(1 + h^{1/2}\log^2 (h/Y)q^{-1}(h,q)^{\frac{1}{2}}\tau(hq)\right).
\]
\end{claim}
\begin{proof}
We analyze $Q_P(k)$ using the Petersson trace formula.
\begin{theorem}[cf.~\cite{IWATOPICS}, page 133, Theorem 9.6]
With all notations as above, let $\mathcal{S}(k)$ be an orthonormal basis for the space of weight-k (holomorphic) modular forms with respect to the Petersson inner product. For $f\in \mathcal{S}(k)$, denote by $a_f(h)$ the $h$'th Fourier coefficient of $f$. Then
\[
\dfrac{\Gamma(k-1)}{(4\pi h)^{k-1}}\sum_{f\in \mathcal{S}(k)}|a_f(h)|^2 = 1 + 2\pi i^{-k}\sum_{c=0(q)}\dfrac{1}{c}S(h,h;c)J_{k-1}(\dfrac{4\pi h}{c}).
\]
\end{theorem}
Although the sum $Q_P(k)$ looks a lot like the spectral side of the Petersson trace formula, there is some difference: the $\lambda_{\phi}(h)$ are the $h$'th Fourier-Whittaker coefficients of an orthonormal basis of the space of min-weight $k$ cusp forms, whereas the $a_f(h)$ are the $h$'th Fourier coefficients of an orthonormal basis of the space of weight $k$ holomorphic modular forms. However, as these two Hilbert spaces are isometric via the isometry (see~\cite{BUMP})
\[
\iota: g(z)\in\text{Span}_{\mathbb{C}}\{f:f\in \mathcal{S}(k)\} \mapsto \dfrac{1}{2\pi}g(x + iy)y^{k/2}e^{in\theta}\in \text{Span}_{\mathbb{C}}\{\phi:\phi\in \mathcal{F}_k\}.
\]
Fixing an orthonormal basis $\mathcal{S}(k)$ for the space of weight $k$ (holomorphic) modular forms, we may assume that our basis $\mathcal{F}_k$ for the space of min-weight $k$ (holomorphic) cusp forms is obtained as the image of $\iota$ on $\mathcal{S}(k)$. A trivial computation shows that the $h$'th Fourier coefficient of $f\in \mathcal{S}(k)$, denoted $a_f(h)$, and the $h$'th Fourier-Whittaker coefficient of $\iota(f)$, denoted $\lambda_f(h)$, are related via:
\[
a_f(h) = (4\pi h)^{(k-1)/2}\lambda_f(h).
\]
Plugging back, we get:
\[
Q_P(k) = (k-1)\left(1 + \sum_{c=0(q)}\dfrac{1}{c}S(h,h;c)J_{k-1}(\dfrac{4\pi h}{c})\right).
\]
We then have (cf.~\cite[14.61]{IWAKOW}):
\[
Q_P(k) \ll k\left(1 + h^{1/2}\log^2 (h/Y)q^{-1}(h,q)^{\frac{1}{2}}\tau(hq)\right),
\]
which completes the proof.
\end{proof}

\subsubsection{Combining the bounds on $R_E(P)$ and $R_D(P)$}
Combining our bounds on $R_D(P)$ and $R_E(P)$ proves that:
\[
R(P) \ll_{\delta} hZ^{-19}Y^{-1-2\delta}\left(1 + h^{\frac{1}{2}}\log^2 (h/Y)\left(Y + Y^{-1}\right)^{\frac{1}{2}}q^{-1}(h,q)^{\frac{1}{2}}\tau(hq)\right),
\]
which completes the proof of Theorem~\ref{theorem:bound_on_R}. 

Theorem~\ref{theorem:bound_on_R} and Theorem~\ref{theorem:bound_on_K} imply that for all $\tau \in SL_2(\mathbb{Z})$:
\[
P(\tau) \ll_{\delta} (\dfrac{h}{Y^{1 + 2\delta}Z^{19}})^{\frac{1}{2}}[1 + h^{\frac{1}{2}}\left(Y + Y^{-1}\right)^{\frac{1}{2}}q^{-1}(h,q)^{\frac{1}{2}}\tau(hq)]^{\frac{1}{2}}\log h/Y,
\]
which completes the proof of Theorem~\ref{theorem:point-wise_bound}.

\section{Equidistribution to prime congruences}
\label{section:equidistribution_to_primes}
In this chapter we prove our main theorem. The key inputs that go into the proof are power saving bounds on linear and bilinear sums. The former was already introduced (see Definition~\ref{definition:linear_sums}), and bilinear sums will be introduced later in this chapter, see Definition~\ref{definition:bilinear_sums}. The corresponding bounds are then propagated into Duke, Friedlander and Iwaniec's sieve, which we use as a black box, see Theorem~\ref{theorem:reformulation}. As a result, we prove our main theorem.

\subsection{Linear sums}
We begin by recalling the definition of linear sums. To ease notations, denote by $\rho_h(n)$ the Weyl sum 
\[
\rho_h(n) := \sum_{\nu\mod n: \nu^2 + 1 = 0(n)}'e(\frac{h\nu}{n}),
\]
where $\sum'$ signifies that we only consider the subset of these $\nu\mod n$ for which the primitive lattice point corresponding to the (root,modulus) pair $(\nu,n)$ lies in our predefined sector $(\alpha,\beta)\subseteq (0,\pi/2)$.

The linear sums, $\mathcal{L}_{d,h}(N)$, are given by:
\[
\mathcal{L}_{d,h}(N) := \sum_{nd\le N}\rho_h(nd).
\]
Our goal in this section is to prove:
\begin{theorem}
\label{theorem:linear_sum}
Let $0 < \delta < 1/4$, one has:
\[
L_{d,h}(N) \ll_{\delta} 
(h,d)^{\frac{1}{52}}\left(\dfrac{N}{d}\right)^{1 + \delta}\left(\dfrac{d^2}{N}\right)^{1/90}.
\]
\end{theorem}
\begin{proof}
The idea is to estimate our linear sums using a smooth version:
\begin{claim}
\label{claim:smooth_sum_bound_g_y_g_theta}
Let $0 < \delta < 1/4$. Suppose $G_y:\mathbb{R}_+\longrightarrow\mathbb{R}_+$ is supported on $[N,2N]$ and has derivatives $G_y^{(j)} \ll N^{-j}$ for $j\le 12$. Let $G_{\theta}:[0,2\pi]\longrightarrow \mathbb{R}$ be supported on $[\alpha + \pi k/2,\beta + \pi k/2], k\in \{0,1,2,3\}$, and has derivatives $G_{\theta}^{(j)} \ll_j Z^{-j}$, where $N^{-1/3} < Z < \dfrac{\beta-\alpha}{2}$, $G_{\theta}$ satisfies $||G_{\theta}^{(j)}||^2 \ll_j Z^{1-2j}$, and its Fourier coefficients, $g_n$, satisfy the bound $|g_n| \ll e^{-2\sqrt{\pi |n| Z}}$. Let $d\ll N$ and $h < N^{1/3}$, then
\[
\bigg|\sum_{\substack{a^2 + b^2 = 0(d)\\ a,b\ge 0\\ (a,b) = 1}}e(\dfrac{h\overline{a}b}{a^2 + b^2})G_y(a^2 + b^2)G_{\theta}(\arctan(\frac{b}{a}))\bigg|  \ll_{\delta} \tau(d)\left(\dfrac{N^{1 + 2\delta}}{Z^{19}}\right)^{\frac{1}{2}}\left[1 + N^{\frac{1}{2}}d^{-1}(h,d)^{\frac{1}{2}}\tau(hd)\right]^{\frac{1}{2}}\log N.
\]
\end{claim}
\begin{proof}
Take $G_y$ and $G_{\theta}$ to be our functions $F$ and $G$, see Definition~\ref{definition:function_F} and Definition~\ref{definition:function_G}. Interpret the left hand side as the appropriate sum of point-wise evaluations of our Poincare series and apply Theorem~\ref{theorem:point-wise_bound} to bound each one of these point-wise valuations.
\end{proof}
Under the further assumption that $d \ll \sqrt{N}$, we may bound:
\[
\bigg|\sum_{\substack{a^2 + b^2 = 0(d)\\ a,b\ge 0\\ (a,b) = 1}}e(\dfrac{h\overline{a}b}{a^2 + b^2})G_y(a^2 + b^2)G_{\theta}(\arctan(\frac{b}{a}))\bigg|  \ll_{\delta} (h,d)^{\frac{1}{4}}\left(\dfrac{N}{d}\right)^{1 + \delta}\left(\dfrac{d^2}{NZ^{38}}\right)^{1/4}
\]
Our goal is to estimate $\mathcal{L}_{d,h}(N)$ using the above smooth summation. In order to do that, we need to bound their difference. Take $g_y$ supported on $[N,2N]$ satisfying $g_y^{(j)} \ll N^{-j}$ for $j\le 12$, and $g_y(n) = \Delta^{12}$ for $N + \Delta N\le n\le 2N - \Delta N$, while $0\le g_y(n) \le \Delta^{12}$ for $0\le n\le \Delta n$ or $2N - \Delta N \le n\le 2N$. Also, take $g_{\theta}$ supported on $[\alpha,\beta]\mod \pi/2$, satisfying $g_{\theta}\left(\arctan(\frac{b}{a})\right) = 1$ for $\alpha + Z \le \arctan(\frac{b}{a}) \le \beta - Z$, $0\le g_{\theta}(a,b)\le 1$ for $\alpha\le \arctan(\frac{b}{a})\le \alpha + Z$, and $\beta - Z\le \arctan(\frac{b}{a}) \le \beta$. Furthermore, assume $g_{\theta}^{(j)} \ll_j Z^{-j}$ and $||g_{\theta}^{(j)}||^2 \ll_j Z^{1-2j}$. Denote the above sum by $\mathcal{L}_{d,h}^*(N)$, then
\[
\mathcal{L}_{d,h}(N) = \dfrac{1}{\Delta^{12}}\mathcal{L}_{d,h}^*(N) + O(|\mathcal{S}(N)|),
\]
where $\mathcal{S}(N)$ is the set of boundary points, defined by
\[
\mathcal{S}(N) := 
\bigg\{(a,b)\in \mathbb{Z}_+\times\mathbb{Z}_+: \left((a,b) = 1\right),\left(a^2 + b^2 = 0(d)\right),
\]
\[
\big(N \le a^2 + b^2 \le N + \Delta N \text{ or } 2N - \Delta N \le a^2 + b^2 \le 
2N\big)\text{ or } \big(\alpha \le \arctan(\dfrac{b}{a})\le \alpha + Z\text{ or } \beta - Z \le \arctan(\dfrac{b}{a})\le \beta\big)\bigg\}.
\]
Before we optimize our choice for $\Delta$, we construct a bound on the cardinality of $\mathcal{S}(N)$.
\begin{claim}
For all $N\in\mathbb{N}$, one has
\[
|\mathcal{S}(N)| \ll_{\delta} d^{\delta}\sqrt{N} + \dfrac{NZ}{d^{1-\delta}} + \dfrac{\Delta N^{1 + \delta}}{d}.
\]
\end{claim}

\begin{proof}
We ``cover" the set $\mathcal{S}(N)$ by:
\[
\mathcal{S}_{\text{radial}}(N) = \Bigg\{\dfrac{\nu}{n}: n=0(d),\left(\left(N\le n\le N + \Delta N\right)\text{ or } \left(2N - \Delta N\le n\le 2N\right)\right), \left(0\le \nu < n\right), \nu^2 + 1 = 0(n)\Bigg\},
\]
and
\[
\mathcal{S}_{\text{angular}}(N) = 
\bigg\{(a,b)\in\mathbb{Z}_+\times \mathbb{Z}_+: N \le a^2 + b^2 \le 2N, a^2 + b^2 = 0(d),
\left(\alpha\le \arctan(\dfrac{b}{a})\le \alpha + Z \text{ or } \beta - Z\le \arctan(\dfrac{b}{a})\le \beta\right)\bigg\}.
\]
Since for each $n\in\mathbb{N}$, $\#\{\nu\mod n: \nu^2 + 1 = 0(n)\} = O(\tau(n)) = O_{\delta}(n^{\delta})$, we bound
\[
|\mathcal{S}_{\text{radial}}(N)| \ll_{\delta} \dfrac{\Delta N^{1 + \delta}}{d}.
\]
As for $S_{\text{angular}}(N)$, we partition this set into two:
\[
\mathcal{S}_{\text{angular}}^{\alpha}(N) = \{(a,b)\in\mathbb{Z}_+\times \mathbb{Z}_+: N \le a^2 + b^2 \le 2N, a^2 + b^2 = 0(d), \alpha\le \arctan(\dfrac{b}{a})\le \alpha + Z\},
\]
and
\[
\mathcal{S}_{\text{angular}}^{\beta}(N) = \{(a,b)\in\mathbb{Z}_+\times \mathbb{Z}_+: N \le a^2 + b^2 \le 2N, a^2 + b^2 = 0(d), \beta - Z\le \arctan(\dfrac{b}{a})\le \beta\}.
\]
Let us shift our focus to the set $\mathcal{S}_{\text{angular}}^{\alpha}(N)$, the analysis for the set $\mathcal{S}_{\text{angular}}^{\beta}(N)$ is essentially the same. We begin by bounding the sector of the annulus prescribed by the angles $[\alpha, \alpha + Z]$ and the radii $\sqrt{N}$ and $\sqrt{2N}$ with the trapezoid $\mathcal{T}^{\alpha}(N)$, bounded by the lines $y = x\tan\alpha$, $y=x\tan(Z + \alpha)$, $x = \sqrt{N}\cos(Z + \alpha)$ and $x = \sqrt{2N}\cos(\alpha)$ (see Figure~\ref{figure:areas}).
\begin{figure}[htp]
    \label{figure:areas}
    \centering
    \includegraphics[height=6cm]{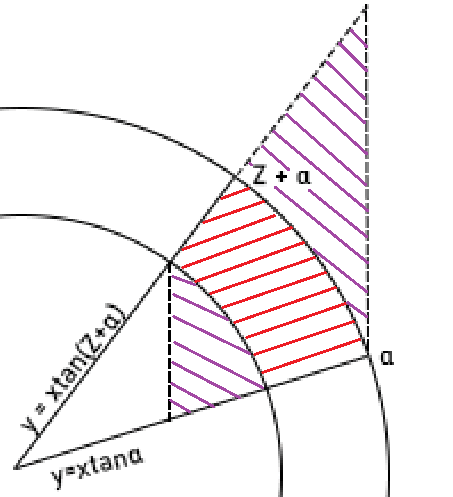}
    \caption{An illustration of the completion of the subsector of the annulus (red area) into the trapezoid (purple and red area).}
    \label{fig:primitive}
\end{figure}
Let 
\[
\mathbb{Z}^2_{\text{mod } d} := \{(a,b)\in\mathbb{Z}^2: a^2 + b^2 = 0(d)\},
\]
it is clear that:
\[
\mathcal{S}_{\text{angular}}^{\alpha}(N)\subseteq
\mathcal{T}^{\alpha}(N)\cap \mathbb{Z}^2_{\text{mod } d}.
\]
The set $\mathcal{T}^{\alpha}(N)\cap \mathbb{Z}^2_{\text{mod } d}$ is easy to analyze in vertical lines. Denote
\[
\xi(n) := \#\{(a,b)\in \mathcal{T}^{\alpha}(N)\cap \mathbb{Z}^2_{\text{mod } d}: a = n\},
\]
then
\[
\#\left(\mathcal{T}^{\alpha}(N)\cap \mathbb{Z}^2_{\text{mod } d}\right) = \sum_{\sqrt{N}\cos(Z + \alpha)\le n\le \sqrt{2N}\cos(\alpha)}\xi(n).
\]
Our claim reduces to:
\begin{claim}
For each integer $n\in \mathbb{Z}\cap [\sqrt{N}\cos(\alpha + Z) ,\sqrt{2N}\cos(\alpha)]$, one has
\[
\xi(n) \ll_{\delta} \dfrac{nZ}{d^{1-\delta}} + d^{\delta}.
\]
\end{claim}

\begin{proof}
Let $n$ be as above. The number of integers satisfying $n\tan(\alpha) \le m \le n\tan(\alpha + Z)$ is bounded above by $n(\tan(\alpha + Z) - \tan(\alpha)) + 1$. Among these integers, the ones which solve
\[
m^2 + n^2 = 0(d),
\]
lie in one of $O_{\delta}(d^{\delta})$ residue classes modulo $d$. Hence,
\[
\xi(n) \ll_{\delta} d^{\delta}\left(\dfrac{n(\tan(\alpha + Z) - \tan(\alpha)) + 1}{d} + 1\right).
\]
Since $Z$ is bounded from $\pi/2 - \alpha$, $\tan(\alpha + Z) - \tan(\alpha) \ll Z$, which completes the proof.
\end{proof}
Combining this with our bound for $|\mathcal{S}_{\text{radial}}(N)|$ completes the proof.
\end{proof}
We now complete the proof of Theorem~\ref{theorem:linear_sum}. Plugging our bound on $O(|\mathcal{S}(N)|)$, we obtain
\[
\mathcal{L}_{d,h}(N) \ll_{\delta} \dfrac{1}{\Delta^{12}}\mathcal{L}_{d,h}^*(N) + d^{\delta}\sqrt{N} + \dfrac{NZ}{d^{1-\delta}} + \dfrac{\Delta N^{1 + \delta}}{d}.
\]
Fixing $\Delta = \left(\dfrac{\mathcal{L}_{d,h}^{*}(N)d}{N}\right)^{1/13}$ and applying the bound from Claim~\ref{claim:smooth_sum_bound_g_y_g_theta} on $\mathcal{L}^*_{d,h}(N)$, we get
\[
|\mathcal{L}_{d,h}(N)| \ll_{\delta} 
(h,d)^{\frac{1}{52}}\left(\dfrac{N}{d}\right)^{1 + \delta}\left(\dfrac{d^2}{NZ^{38}}\right)^{1/52}
+ \dfrac{NZ}{d^{1-\delta}} + d^{\delta}\sqrt{N}.
\]
Fixing $Z = (d^2/N)^{1/90}$, we conclude that
\[
|\mathcal{L}_{d,h}(N)| \ll_{\delta} 
(h,d)^{\frac{1}{52}}\left(\dfrac{N}{d}\right)^{1 + \delta}\left(\dfrac{d^2}{N}\right)^{1/90},
\]
which completes the proof.
\end{proof}

\subsection{Bilinear sums}
\begin{definition}
\label{definition:bilinear_sums}
Let $(\rho_h(n))_{n=1}^{\infty}$ be our sequence of Weyl sums described above. Let $M, P$ be natural numbers, and let $(\alpha_m)_{m=1}^{\infty}$ and $(\beta_p)_{p=1}^{\infty}$ be sequences of complex numbers, where the sequence $(\beta_p)_{p=1}^{\infty}$ is supported only on prime indices. We define the bilinear sum $B(M,P)$ of the sequence $(\rho_h(n))_{n=1}^{\infty}$ by
\[
B(M,P) := \sum_{M\le m\le 2M}\sum_{P\le p\le 2P}\alpha_m\beta_px_{mn}.
\]
\end{definition}
Our goal in this section is to prove:
\begin{claim}
\label{claim:bilinear_sums}
Fix $h \in \mathbb{Z}$ and $0 < \epsilon < 1/3$. Let $M, P, N\in\mathbb{N}$ satisfy $M\ge N^{2/3}, P\le N^{1/3-\epsilon}$. Assume both sequences $(\alpha_m)_{m=1}^{\infty}$ and $(\beta_n)_{n=1}^{\infty}$ are supported in the intervals $[M,2M]$ and $[P,2P]$, respectively, and that $(\beta_n)_{n=1}^{\infty}$ is supported on primes. Then, for all $0 < \delta < 1/4$ and $N$ large enough, one has
\[
|B(M,P)| \ll_{\delta} ||\alpha||\cdot||\beta||\cdot \left(M^{1/2} + M^{41/84 + \delta}P^{11/21}\right).
\]
\end{claim}
\begin{proof}
We would like to replicate~\cite[Proposition 2]{DFI} verbatim, however, there are two problems. The first is concerned with our admissible range for $h$; we explain. In~\cite[Proposition 2]{DFI}, bounding the bilinear sum involves an estimation of various smooth linear sums. Since the admissible range for $h$ in Duke, Friedlander and Iwaniec's bound on (smooth) linear sums is $h \ll N$, we must make sure that we only apply Claim~\ref{claim:smooth_sum_bound_g_y_g_theta} whenever $h < N^{1/3}$.

Duke, Friedlander and Iwaniec invoke their bound on smooth linear sums with $N \gg MP^2$ and with parameter ``$h$" equal to $hP$. To make sure that we can apply our bound whenever they do, we must make sure that
\[
hP < (MP^2)^{1/3}\iff h < \left(\dfrac{M}{P}\right)^{1/3},
\]
which follows from $h \ll M^{1/6}$. However, since $h$ is a constant and $M$ tends to $\infty$ with $N$, the assumptions of Claim~\ref{claim:smooth_sum_bound_g_y_g_theta} hold whenever $N$ is large enough.

The second issue is that in order to adapt the proof of Proposition 2 in~\cite{DFI} to our setting, we do not require bounds on smooth sums precisely, but rather on ``half smooth" sums, which are sums of the form:
\[
A_{p_1,p_2}(M) = \sum_{m = M/2}^{4M}g(m)\rho_{h(p_2-p_1)}(mp_1p_2),
\]
where $g(m)$ is any smooth, nonnegative function, which majorizes the indicator of the interval $[M,2M]$. If we take $g$ to be as in Claim~\ref{claim:smooth_sum_bound_g_y_g_theta}, then we may bound $A_{p_1,p_2}(M)$ using
$|A_{p_1,p_2}(M)| \ll |\mathcal{L}_{d,h}^*(MP^2)| + O(|\mathcal{S}_{\text{angular}}(MP^2)|)$. Following the ideas in the proof of Theorem~\ref{theorem:linear_sum}, we find that for all $0 < \delta < 1/4$:
\[
|A_{p_1,p_2}(M)| \ll_{\delta} M^{1 + \delta}\left(\dfrac{P^2}{MZ^{38}}\right)^{1/4} + P^{\delta}\sqrt{MP^2} + MP^{\delta}Z.
\]
Fixing $Z = \left(\dfrac{P^2}{M}\right)^{1/42}$, yields
\[
|A_{p_1p_2}(M)| \ll_{\delta}M^{1 + \delta}\left(\dfrac{P^2}{M}\right)^{1/42}.
\]
Equipped with this bound, following the rest of the proof of Proposition 2 in~\cite{DFI} gives
\[
|B(M,P)| \ll_{\delta} ||\alpha||\cdot||\beta||\cdot \left(M^{1/2} + M^{41/84 + \delta}P^{11/21}\right).
\]
which is what we wanted to show.
\end{proof}

\subsection{Duke, Friedlander and Iwaniec's sieve}
We pack Duke, Friedlander and Iwaniec's sieving argument into a black box (cf. \cite[p.433-437]{DFI}).
\begin{theorem}
\label{theorem:reformulation}
Assume that for all $h\in \mathbb{Z}$, $N\in\mathbb{N}$ large enough and $0 < \epsilon < 1/3$:
\begin{itemize}
    \item (A.) One has
    \[
    \left|\sum_{d \le N^{1/2 - \epsilon}}\lambda_d\mathcal{L}_{d,h}(N)\right| \ll_{\epsilon} \dfrac{N}{\log^2 N}\max_{d\le N}\{|\lambda_d|\}.
    \]
    \item (B.) For all integers $P,Q$ satisfying $N^{(\log\log N)^{-3}} < Q < P < N^{1/3 - \epsilon}$, 
    and any sequences of complex numbers $(\alpha_m)_{m=1}^{\infty}$, and $(\beta_n)_{n=1}^{\infty}$, with the latter supported on prime indices, such that $|\alpha_m|\le \omega(m)$ and $|\beta_n| \le 1$, one has
    \[
    \left|\sum_{Q\le n\le P}\beta_n\sum_{m\le N/n}\alpha_m\rho_h(mn)\right| \ll_{\epsilon} \dfrac{N}{\log^{10}N}.
    \]
\end{itemize}
Then
\[
\sum_{p\le N}\rho_h(p) \le \epsilon\pi(N).
\]
\end{theorem}
Our main theorem reduces to showing that the sequence $(\rho_h(n))_{n=1}^{\infty}$ has properties $(A.)$ and $(B.)$.

We prove a stronger version of property $(A.)$:
\begin{claim}
There exists an $\eta > 0$, such that
\[
\sum_{d\le N^{1/2 - \epsilon}}|\mathcal{L}_{d,h}(N)| \ll_{\epsilon} N^{1-\eta}.
\]
\end{claim}
\begin{proof}
Let $0 < \delta < 1/4$ be decided upon in hindsight, by Theorem~\ref{theorem:linear_sum}, we have
\[
|\mathcal{L}_{d,h}(N)| \ll_{\delta} \left(\dfrac{N}{d}\right)^{1 + \delta}\left(\dfrac{d^2}{N}\right)^{1/90}.
\]
Plugging back, we get
\[
\sum_{d\le N^{1/2 - \epsilon}}|\mathcal{L}_{d,h}(N)| \ll_{\delta} \sum_{d\le N^{1/2 - \epsilon}}\left(\dfrac{N}{d}\right)^{1 + \delta}\left(\dfrac{d^2}{N}\right)^{1/90} \ll_{\delta} N^{1 + \delta - \epsilon/45}.
\]
Therefore, setting $\delta = \epsilon/90$, our claim holds for $\eta = \epsilon/90$.
\end{proof}

Next, we prove: 
\begin{claim}
Property $(B.)$ holds.
\end{claim}
\begin{proof}
All notations as above. An examination of the proof of Claim~\ref{claim:bilinear_sums} (covering the hyperbola with rectangles) shows that for all $0\le \delta < 1/4$:
\[
\left|\sum_{2^{-k-1}P\le n\le 2^{-k}P}\beta_n\sum_{m\le N/n}\alpha_m\rho_h(mn)\right| 
\ll_{\delta,A} ||\alpha||\cdot||\beta||\cdot \left((2^kN\log N/P)^{1/2} + (2^kN/P)^{41/84 + \delta}(P/2^k)^{11/21} + (N/\log^A N)\right).
\]
Since for all $m > 1$: $\omega(m) \ll \log m$, we have
\[
||\alpha|| \ll 2^{k/2}\sqrt{N/P}\log N,\quad ||\beta|| \ll 2^{-k/2}\sqrt{P}.
\]
Hence,
\[
\left|\sum_{2^{-k-1}P\le n\le 2^{-k}P}\beta_n\sum_{m\le N/n}\alpha_m\rho_h(mn)\right| 
\ll_{\delta, A} 2^{k/2}N\log^{3/2}N/P^{1/2} + N^{83/84 + \delta}P^{3/84}/2^{11k/21} + (N/\log^A N).
\]
Summing from $k = 0$ to $\lceil\log_2(P/Q)\rceil - 1$, we obtain the bound
\[
\left|\sum_{Q\le n\le P}\beta_n\sum_{m\le N/n}\alpha_m\rho_h(mn)\right| \ll_{\delta, A} (N/\log^{A-1} N) + \sum_{k=0}^{\lceil\log_2(P/Q)\rceil - 1}2^{k/2}N\log^{3/2}N/P^{1/2} + N^{83/84 + \delta}P^{3/84}/2^{11k/21} 
\]
\[
\ll_{\delta, A} N/\log^{A-1} N + N\log^{3/2}N/Q^{1/2} + N^{83/84 + \delta}P^{3/84}.
\]
Since $Q > N^{(\log\log N)^{-3}}$ and $P < N^{1/3 - \epsilon}$, by adjusting $A = 11$ and $\delta < 3\epsilon/84$, we obtain:
\[
\left|\sum_{Q\le n\le P}\beta_n\sum_{m\le N/n}\alpha_m\rho_h(mn)\right| \ll_{\epsilon} N/\log^{10} N,
\]
which completes the proof.
\end{proof}
Since both properties (A.) and (B.) hold, Theorem~\ref{theorem:reformulation} implies that $\sum_{p\le N}\rho_h(p) \le \epsilon\pi(N)$, which completes the proof of our main theorem.

\section{Constructing $\rho(t,n)$}
\label{appendix:constructing_rho}
We begin by stating the main theorem concerning $\rho(t,n)$.
\begin{theorem}
\label{claim:rho_exists}
There exists a function $\rho(t,n):\mathbb{C}\times\mathbb{Z}\longrightarrow\mathbb{C}$ satisfying the following properties.
\begin{itemize}
    \item (1) $\rho(t,n) > 0$ for all $(t,n)$ pairs of type and weight coming from eigenforms of the Laplacian on $\Gamma_0(N)\char`\\SL_2(\mathbb{R})$ with $N\in\mathbb{N}$.
    \item (2) $\rho(t,n)$ is the Harish-Chandra/Selberg transform of a point-pair invariant $k(g,h)$ of the form $\sum_{n\in 2\mathbb{Z}}k_n(g,h)$, satisfying the absolute convergence property.
    \item (3) (Real valued): $k(g,h)$ attains real values.
    \item (4) (Positivity): For all $g,h\in SL_2(\mathbb{R})$ one has $k_0(g,h) \ge \sum_{n\neq 0}|k_n(g,h)|$.
    \item (5) $\rho(t,n)^{-1}$ is bounded by a polynomial in $t,n$.
\end{itemize}
\end{theorem}

The goal of this section is to give a brief survey of the construction of $\rho(t,n)$, which played a crucial role in this work. See~\cite{SAARTHESIS} for a full account of this construction.

The key difficulty in showing that a specific choice of $\rho$ satisfies conditions $(1)-(5)$ is that some of these properties, such as $(3)$ and $(4)$, for example, do not concern $\rho$ directly, but rather its inverse Harish-Chandra/Selberg transform. Moreover, it is not immediate that an arbitrary choice of $\rho$ should even be in the image of the Harish-Chandra/Selberg transform (hence property $(2)$). In fact, showing that our particular choice of $\rho$ is in the image of the Harish-Chandra/Selberg transform was one of the main challenges in the construction.

Before we introduce our explicit choice of $\rho$, we motivate our choices by recalling the definition of the inverse Harish-Chandra/Selberg transform. There isn't a non-spherical inverse Selberg transform per se, but instead, there is one for each weight, and the spherical inverse Selberg transform is just the weight zero case of this family of operators.

\subsection{Inverting the Harish-Chandra/Selberg transform}
\label{section:inversion}
We first recall a few definitions given in Section~\ref{section:preliminary_discussion_k}.

\begin{definition}[\cite{HEJHAL1} page 357, Definition 2.3]
Let $z\in\mathbb{H}$, and $\sigma\in SL_2(\mathbb{R})$ given by
\[
\sigma =
\begin{pmatrix}
a& b\\
c& d
\end{pmatrix}.
\]
We define
$j_{\sigma}(z) := \dfrac{cz + d}{|cz + d|} = e^{i\arg (cz + d)}$.
\end{definition}

\begin{definition}[\cite{HEJHAL1} page 359, Definition 2.9]
Let $n\in \mathbb{Z}$, we define
\[
\forall z,w\in\mathbb{H}:\quad H_n(z,w) := i^n\dfrac{(w-\overline{z})^n}{|w-\overline{z}|^n}.
\]
\end{definition}

\begin{definition}[cf.~\cite{HEJHAL1} page 359, Definition 2.10)]
Let $\Phi_n:\mathbb{R}_+\longrightarrow\mathbb{R}$ be any $C^2$-class function. We define the weight-$n$ point-pair invariant attached to $\Phi_n$ to be
\[
k_n(g,h) = \Phi_n(u(\pi(g),\pi(h)))H_n(\pi(g),\pi(h))e^{in(\theta_g - \theta_h)},
\]
where $u$ is the hyperbolic distance on $\mathbb{H}$.
\end{definition}

We also require the definition of the weight-$n$ Harish-Chandra/Selberg transform, (cf. Definition~\ref{definition:harish-chandra_selberg_transform}).

\begin{definition}
\label{definition:selberg_transform}
Given a weight-$n$ point-pair invariant $k_n(g,h)$, we define its weight-$n$ Harish-Chandra/Selberg transform $\rho_{k_n}(t)$ by
\[
\rho_{k_n}(t) = \int_{SL_2(\mathbb{R})}k_n(I,h)y(h)^{1/2 + it}e^{in\theta(h)}\dfrac{dxdyd\theta}{y(h)^2},
\]
where $y(h)$ and $\theta(h)$ are the Iwasawa $y$ and $\theta$ coordinates of $h$, respectively.
\end{definition}

Under certain conditions on the weight-$n$ point-pair invariant $k_n(g,h)$, the weight-$n$ Harish-Chandra/Selberg transform is invertible. We name this subset of functions ``admissible weight-$n$ point-pair invariants".
\begin{definition}[\cite{HEJHAL1}]
\label{definition:weight_m_point_pair}
Let $k_n(g,h)$ be a weight-$n$ point-pair invariant. We say that it is an admissible weight-$n$ point-pair invariant, if the associated $\Phi_n$ is in $C^4(\mathbb{R}_+)$, and
\[
|\Phi_n^{(k)}(t)| \ll (t + 1)^{-\tau - k},
\]
where $\tau = \max\{|n|/2,1\}$, $k=0,1,2,3,4$.
\end{definition}

Following Hejhal (see~\cite[p. 386]{HEJHAL2}), the weight-$n$ Harish-Chandra/Selberg transform $\rho_{k}(t)$ of $k_m(g,h)$ is computed using the following 3-step recipe.

Fix $Q_n:\mathbb{R}_+\longrightarrow\mathbb{C}$ by
\[
\forall w \ge 0:\quad Q_n(w) = \int_{-\infty}^{\infty}\Phi_n(w + v^2)\cdot\left(\dfrac{\sqrt{w + 4} + iv}{\sqrt{w + 4} - iv}\right)^{n/2}dv.
\]
When $\Phi_n$ comes from an admissible weight-$n$ point-pair invariant, $\Phi_n$ may be recovered from $Q_n$ by
\[
\forall x\ge 0:\quad \Phi_n(x) = -\dfrac{1}{\pi}\int_{-\infty}^{\infty}Q_n'(x + r^2)\cdot\left(\dfrac{\sqrt{x + 4 + r^2} - r}{\sqrt{x + 4 + r^2} + r}\right)^{n/2}dr,
\]
where $Q_n'$ is the derivative of $Q_n$. Let $g_n(u)$ be defined by
\[
g_n(u) = Q_n(e^u + e^{-u} - 2).
\]
Alternatively,
\[
Q_n(x) = g_n(\text{arccosh}(x/2 + 1)).
\]
The Selberg transform of $k_n(g,h)$ is given by the Fourier transform of $g_n(u)$.
\[
\rho_{k_n}(t) = \int_{-\infty}^{\infty}g_n(u)e^{itu}du.
\]
By Fourier inversion, for $\rho_{k_n}(t)$ of $C^{1}(\mathbb{R})$-class, we also have
\[
g_n(u) = \int_{-\infty}^{\infty}\rho_{k_n}(t)e^{-itu}dt.
\]

\subsection{Fixing our choice of $\rho$}
In the previous section we learned that a sufficient condition for a function $\rho(t,n)$ to be in the image of the Harish-Chandra/Selberg transform, is that its inverse Selberg transform has an associated function $\Phi(x)$, which is $C^4$ and satisfies certain growth conditions.

In~\cite{DFI}, the authors construct a test function $h(t)$ (recall that in the spherical case $n = 0$) of the form $(1 + t^2)^{-1} - (4 + t^2)^{-1}$. This function has the miraculous property of yielding a positive point-pair invariant, and as we are also interested in this property in our non-spherical setting, we might try and naively consider $h(t)$ or a minor variant thereof, as a first choice for our test function. This choice has several issues, and explaining what these are is instructive. 

Let's examine the simplest issue concerning this choice of $\rho$. Our point-pair invariant, $k(g,h)$, which we defined as a sum of the inverse Selberg transforms of $\rho(t,n)$ for every even integer $n$, may not converge. These issues could be easily fixed by adjusting $\rho(t,n)$ so that it exhibits polynomial decay in $n$. This may lead us to examine a family of test functions of the form
\[
\rho(t,n) = ((X_n)^2 + t^2)^{-1} - (4(X_n)^2 + t^2)^{-1},
\]
where $X_n$ is a sequence of real numbers which grows fast enough to ensure the convergence of $k(g,h)$ and $K(g)$. Fixing $X_n \gg |n|$ would be enough for this task. However, it so happens that the associated $\Phi$ turns out to be $C^0$ instead of $C^4$-class. To understand why this comes to be, we reexamine the Selberg inversion operator.

Starting with a function $\rho(t,n)$, we denote by $g_n(x)$ its Fourier inversion with respect to the variable $t$, and denote by $Q_n(x)$ the function defined by $Q_n(x) := g_n(\text{arccosh}(x/2 + 1))$. The function $\Phi(x)$ is derived from $Q_n(x)$ via
\[
\Phi_n(x) = -\dfrac{1}{\pi}\int_{-\infty}^{\infty}Q_n'(x + r^2)\cdot\left(\dfrac{\sqrt{x + 4 + r^2} - r}{\sqrt{x + 4 + r^2} + r}\right)^{n/2}dr.
\]
A sufficient condition for $\Phi_n(x)$ to be $C^4$ is that $Q_n(x)$ is $C^5$, and that the integral converges. To see why we can only show that $\Phi(x)$ is $C^0$, we will show that $Q_n(x)$ is only $C^1$.

Let us consider the family of functions
\[
\rho_{f}(t) = \dfrac{1}{f^2 + t^2}.
\]
$\rho_f(t)$ has the Fourier transform
\[
g(u) = \int_{-\infty}^{\infty}\dfrac{e^{-itu}dt}{f^2 + t^2} = \dfrac{\pi}{f}e^{-f|u|}.
\]
Denoting by $Q_f(x)$ the associated function $Q$, we have
\[
Q_f(x) = g(\text{arccosh}(x/2 + 1)) = \dfrac{\pi}{f}e^{-f\cdot\text{arccosh}(x/2 + 1)} = \dfrac{\pi}{f}\cdot\dfrac{1}{\left(1 + \dfrac{x}{2} + \dfrac{\sqrt{x}\sqrt{x + 4}}{2}\right)^f}.
\]
$Q_f(x)$ has the following Puiseux series:
\[
\dfrac{Q_f(x)}{\pi} = \dfrac{1}{f} - \sqrt{x} + \dfrac{fx}{2} + \dfrac{1}{24}(1 - 4f^2)x^{3/2} + O\left((f^3 + 1)x^2\right) \dfrac{1}{24}f(f^2 - 1)x^2 + O\left((f^4 + 1)x^{5/2}\right),
\]
and by linearity, we see that the coefficient of $x^{3/2}$ of the associated $Q$-function associated to $\rho(t,n) = (X_n^2 + t^2)^{-1} - (4X_n^2 + t^2)^{-1}$ is $\pi(16X_n^2 - 4X_n^2) \neq 0$, which yields an obstruction for the $C^1$-ness of $\Phi(x)$.

We may try to remedy this issue by replacing our previous $\rho$ by 
\[
\rho(t,n) = ((aX_n)^2 + t^2)^{-1} + ((bX_n)^2 + t^2)^{-1} - ((cX_n)^2 + t^2)^{-1} - ((dX_n)^2 + t^2)^{-1},
\]
and fixing $a,b,c,d$ such that $a,b,c,d > 0$, and the (super)sets $\{a,b\} \neq \{c,d\}$ (so that $\rho \not\equiv 0$), such that the equation
\[
a^2 + b^2 = c^2 + d^2
\]
is satisfied.

Our new, 4-parameter $\rho$ has an associated $Q(x)$ which is of $C^2$-class. One may naively try and consider a function of the form $((aX_n)^2 + t^2)^{-1} - ((bX_n)^2 + t^2)^{-1}$ instead, however, in this case, the vanishing of the coefficient of $x^{3/2}$ is equivalent to $a^2 = b^2$, which would render $\rho = 0$.

Similarly, in order for the associated $Q(x)$ to be $C^5$, one obtains analogue conditions on higher (even) powers of $a,b,c,d$. In order to avoid getting an overdetermined system of algebraic equations, and thus risking the associated $\rho$ being zero, we need to introduce more variables. Pushing this idea all the way up to the $C^5$-ness of the associated $Q(x)$ has led us to consider a 1-parameter family of functions.

\begin{definition}
\label{definition:rho}
Let
\[
\rho_X(t) = \dfrac{1}{t^2 + (aX)^2} + \dfrac{1}{t^2 + (bX)^2} + \dfrac{1}{t^2 + (cX)^2} + \dfrac{1}{t^2 + (dX)^2} + \dfrac{1}{t^2 + (2.5X)^2}
\]
\[
- \dfrac{1}{t^2 + X^2} - \dfrac{1}{t^2 + 4X^2} - \dfrac{1}{t^2 + 9X^2} - \dfrac{1}{t^2 + 16X^2} - \dfrac{1}{t^2 + 25X^2},
\]
where $a,b,c,d$ are the 4 positive real solutions, unique (up to permutation), of the system
\[
a^2 + b^2 + c^2 + d^2 + (2.5)^2 = 1^2 + 2^2 + 3^2 + 4^2 + 5^2 = 55,
\]
\[
a^4 + b^4 + c^4 + d^4 + (2.5)^4 = 1^4 + 2^4 + 3^4 + 4^4 + 5^4 = 979,
\]
\[
a^6 + b^6 + c^6 + d^6 + (2.5)^6 = 1^6 + 2^6 + 3^6 + 4^6 + 5^6 = 20515,
\]
\[
a^8 + b^8 + c^8 + d^8 + (2.5)^8 = 1^8 + 2^8 + 3^8 + 4^8 + 5^8 = 462979.
\]
Mathematica gives them as
\[
a^2 \simeq 0.47932,\quad b^2 \simeq 6.87175,\quad c^2\simeq 16.4822,\quad d^2\simeq 24.9163.
\]
\end{definition}
\begin{claim}
The function $Q(x)$ associated to $\rho_X(t)$ is $C^5$.
\end{claim}
See~\cite{SAARTHESIS} for a proof of this fact, as well as a complete proof that this $\rho$ induces an admissible weight $n$ point-pair invariant for large enough $X$.

We ultimately end up choosing $\rho(t,n) = c_n\rho_{X_n}(t)$ for some sequence of positive $c_n$ and $X_n$.
\begin{definition}
\label{definition:final_rho}
For every $t\in\mathbb{C}$ and $n\in\mathbb{Z}$, we define 
\[
\rho(t,n) =
\begin{cases}
\rho_{X_0}(t)& n=0,\\
\dfrac{C}{1200}\dfrac{1}{n^{10} + 1}\rho_{X_n}(t)& n\in 2\mathbb{Z}\setminus\{0\},\\
0& \text{otherwise}.
\end{cases}
\]
Where $\rho_X(t)$ is as in Definition~\ref{definition:rho}, $X_n = (|n| + 2)X_0$, $X_0 = 100$, and $C > 0$ satisfies:
\[
|k_0(g,h)| \ge C|k_m(g,h)|/X_m
\]
for all $m\in 2\mathbb{Z}$ and $g,h\in SL_2(\mathbb{R})$; with $k_n(g,h)$ being the weight-$n$ inverse Selberg transform of $\rho_{X_n}(t)$.
\end{definition}
In first sight, it is totally unclear that such a number $C$ should exist. We refer the reader to the discussion in~\cite{SAARTHESIS}.

Now that we have defined $\rho(t,n)$ we claim:
\begin{claim}
\label{claim:validity}
The function $\rho(t,n)$ satisfies conditions $(1)-(5)$.
\end{claim}
\begin{proof}
See~\cite{SAARTHESIS}.
\end{proof}

\printbibliography
\end{document}